\documentclass{amsart}

\usepackage{amsmath}
\usepackage{amssymb}
\usepackage{amsthm}
\usepackage{mathrsfs}
\usepackage{graphicx}

\title[On Fuchsian Loci of ${\rm PSL}_n(\mathbb{R})$-Hitchin components]{On Fuchsian loci of ${\rm PSL}_n(\mathbb{R})$-Hitchin components of a pair of pants}
\author{Yusuke Inagaki}
\address{Graduate~School~of~Science, Osaka~University}
\email{y-inagaki@cr.math.sci.osaka-u.ac.jp}
\begin{document}
\maketitle

\newtheorem{theorem}{Theorem}[section]
\newtheorem{corollary}[theorem]{Corollary}
\newtheorem{lemma}[theorem]{Lemma}
\newtheorem{proposition}[theorem]{Proposition}
\newtheorem{definition}[theorem]{Definition}
\newtheorem{question}[theorem]{Question}
\newtheorem{example}[theorem]{Example}
\newtheorem{remark}[theorem]{Remark}

\begin{abstract}
In this paper, we study Fuchsian loci of ${\rm PSL}_n(\mathbb{R})$-Hitchin components. In particular, using the Bonahon-Dreyer's parametrization of ${\rm PSL}_n(\mathbb{R})$-Hitchin components, we give an explicit parametrization of Fuchsian loci of a pair of pants.    
\end{abstract}

\section{Introduction}
Let $S$ be a compact oriented surface with negative Euler characteristic number.  
The ${\rm PSL}_n(\mathbb{R})$-{\it Hitchin component} of $S$ is a connected component of the ${\rm PSL}_n(\mathbb{R})$-character variety
\[ \mathcal{X}_n(S) = {\rm Hom}(\pi_1(S), {\rm PSL}_n(\mathbb{R}))//{\rm PSL}_n(\mathbb{R}). \]
Recall that the {\it Teichm\"uller space} of $S$, denoted by $\mathscr{T}(S)$, which is the space of marked hyperbolic structures on $S$, is identified with a connected component of $\mathcal{X}_2(S)$ (see \cite{Go} and \cite{IT}). 
In particular this component consists of discrete faithful representations $\rho$ which induce orientation-preserving isometries $S \cong \mathbb{H}^2/\rho(\pi_1(S))$.
We call this component the {\it Teichm\"uller component}, which we denote it by ${\rm Fuch}_2(S)$, and representations $\rho: \pi_1(S) \rightarrow {\rm PSL}_2(\mathbb{R})$ with $[\rho] \in {\rm Fuch}_2(S)$ {\it Fuchsian representations}. 
Let $\iota_n : {\rm PSL}_2(\mathbb{R}) \rightarrow {\rm PSL}_n(\mathbb{R})$ be the canonical $n$-dimensional irreducible representation. 
This representation induces a map $(\iota_n)_* : \mathcal{X}_2(S) \rightarrow \mathcal{X}_n(S) $ by the left composition. 
The ${\rm PSL}_n(\mathbb{R})$-{\it Hitchin component} of $S$, denoted by ${\rm Hit}_n(S)$, is the connected component of $\mathcal{X}_n(S)$ which contains the image $(\iota_n)_*({\rm Fuch}_2(S))$ of the Teich\"uller component. 
Hitchin components were originally introduced in \cite{Hi} and they are now studied in several contexts such as a generalization of Teichm\"uller theory, and deformation theory of some geometric structures. 
For example, it is known that ${\rm Hit}_n(S)$ is diffeomorphic to a Euclidian space as in the case of Teichm\"uller spaces (\cite{Hi}), and ${\rm Hit}_n(S)$ coincides with the deformation space of convex real projective structures on $S$ when $n=3$ (\cite{CG}), and with the deformation space of convex foliated projective structures on the unit tangent bundle $T^1S$ when $n=4$ (\cite{GW2}).

In this paper we focus on {\it Fuchsian loci} of Hitchin components.
The {\it Fuchsian locus} in ${\rm Hit}_n(S)$ is the subset defined by ${\rm Fuch}_n(S)=(\iota_n)_*({\rm Fuch}_2(S))$.  
It is known that Fuchsian loci have some interesting properties with respect to entropies and metrics on ${\rm Hit}_n(S)$ (see \cite{BCFS}, \cite{LW}, \cite{PS} and \cite{Z}).
However it is unknown that how Fuchsian loci lie in Hitchin components concretely. 
This paper gives an answer of this natural question in the case of a pair of pants by parameterizing Fuchsian loci. 
For the parameterization, we use Bonahon-Dreyer's work \cite{BD2} in which they constructed a homeomorphic parameterization of ${\rm Hit}_n(S)$ by a convex polytope $\mathcal{P}$ in a Euclidian space. 
Under their parameterization, we detect the image of ${\rm Fuch}_n(S)$ for general $n \geq 3$ in the polytope $\mathcal{P}$ in the case of a pair of pants, which is an explicit description of ${\rm Fuch}_n(S)$.

Here we briefly review a history of studies of the Anosov property of elements in ${\rm Hit}_n(S)$, called {\it Hitchin representations}, which is used in the construction of the Bonahon-Dreyer parameterization. 
In \cite{Hi}, Hitchin studied components of ${\rm PSL}_n(\mathbb{R})$-character varieties $\mathcal{X}_n(S)$ and defined Hitchin components.
He however commented ``{\it Unfortunately the analytic point of view used for the proofs gives no indication of the geometric significance of the Teichm\"uller component}'' in the introduction. 
(In \cite{Hi}, Hitchin used the term Teichm\"uller components for Hitchin components.)
After Hitchin's work, Labourie\cite{L} found geometric and dynamical meanings of Hitchin representations.
He studied representations which induce flat associated bundles with dynamical properties and showed that Hitchin representations were such representations. 
These representations are now called {\it Anosov representations} and developed into various directions (see \cite{GGKW}, \cite{GW1} and \cite{KLP}). 
We remark that Fock-Goncharov \cite{FG} also studied properties of Hitchin representations from a viewpoint of cluster theory and found {\it positive properties} of Hitchin representations around the same time. 
The dynamical properties of Anosov representations are often called {\it Anosov properties}. 
An important point of the Anosov properties is that we can take nice curves in flag varieties associated to Anosov representations. 
More concretely, for any Hitchin representation $\rho : \pi_1(S) \rightarrow {\rm PSL}_n(\mathbb{R})$, there exists a unique continuous $\rho$-equivariant map $\xi_{\rho}: \partial_{\infty}\pi_1(S) \rightarrow {\rm Flag}(\mathbb{R}^n)$ with the {\it hyperconvex property} (see \cite{Gu} and \cite{L}) where $\partial_{\infty}\pi_1(S)$ is the Gromov boundary of $\pi_1(S)$ and ${\rm Flag}(\mathbb{R}^n)$ is the set of flags in $\mathbb{R}^n$. 
These curves are called by several names such as {\it flag curves, limit curves} and {\it Anosov maps}. In this paper we call such curves flag curves. 

Bonahon-Dreyer used the Anosov property of Hitchin representations to parametrize ${\rm Hit}_n(S)$. 
Let us explain an outline of the construction of their parametrization and prepare some necessary notions to define it at first.
We endow a marked hyperbolic structure $S$ with $m$. 
An {\it $m$-geodesic lamination} on $S$ is a disjoint family of simple $m$-geodesic curves on $S$ which are either biinfinite curves, called {\it biinfinite leaves}, or closed curves, called {\it closed leaves}. 
We remark that there exists a natural bijection between the set of $m$-geodesic laminations and the set of $m'$-geodesic laminations for different hyperbolic structures $m$ and $m'$ of $S$.
Therefore we call an $m$-geodesic lamination a geodesic lamination simply.
A geodesic lamination is said to be {\it maximal} if it is contained in no other geodesic lamination properly.
In this paper we consider only laminations consisting of finitely many leaves which induce an ideal triangulation on $S$. 
(General cases are discussed in \cite{BD1}.) 
Fix a maximal lamination $\mathcal{L}$ on $S$ whose biinfinite leaves and closed leaves are denoted by $h_1, \cdots, h_s$ and $g_1 \cdots, g_t$ respectively.
The lamination $\mathcal{L}$ induces an ideal triangulation on $S$ which consists of complementary ideal triangles $T_1, \cdots, T_u$.
Under these settings, Bonahon-Dreyer constructed two invariants called {\it shearing invariants} and {\it triangle invariants} of Hitchin representations $\rho$, which are valued in $\mathbb{R}$.
The shearing invariants of $\rho$ are defined for leaves $\{ h_i, g_i \}$, denoted by $\sigma_p^{\rho}(h_i)$ and $\sigma_p^{\rho}(g_i)$ where $p$ is an integer index with $1 \leq p \leq n-1$, and the triangle invariants of $\rho$ are defined for ideal triangles $\{ T_i\}$, denoted by $\tau_{pqr}^{\rho}(T_i, v)$ where $v$ is an ideal vertex of $T_i$ and $pqr$ is an index such that $p,q,r \geq 1$ and $p+q+r=n$. 
Bonahon-Dreyer showed that the map $\Phi_\mathcal{L}$ sending $[\rho] \in {\rm Hit}_n(S)$ to a point in a Euclidian space, which is defined $\tau_{pqr}^{\rho}(T_i, v)$, $\sigma_p^{\rho}(h_i)$ and $\sigma_p^{\rho}(g_i)$ of $\rho$ for all $T_i$, $v$, $h_i$, $g_i$, is an onto-homeomorphism whose image is a convex polytope. 
More details of the construction of the Bonahon-Dreyer parameterization will be given in Section 4.
The goal of this paper is to detect the image of the Fuchsian locus in ${\rm Hit}_n(S)$ by $\Phi_\mathcal{L}$, i.e. to compute triangle and shearing invariants of Hitchin representations in Fuchsian loci of a pair of pants concretely.

The result of our computation is as follows. 
At first, we parametrize the Fuchsian representations of the fundamental group $\pi_1(P)$ of a pair of pants $P$ which are identified with hyperbolic structures of $P$. 
It is well known in hyperbolic geometry that hyperbolic structures on $P$, which make the boundary of $P$ totally geodesic, are uniquely determined by the hyperbolic lengths of the boundary components of $P$ up to isometry.
Therefore we can parameterize the Fuchsian representations by the data of the hyperbolic lengths of the boundary components of $P$.

\begin{proposition}[Proposition 5.1]
Let ${\bf m}=(l_A, l_B, l_C)$ be a triple of the hyperbolic lengths of the boundary components $A,B$ and $C$ of $P$. 
Then Fuchsian representations associated to ${\bf m}$ are conjugate to the representation given as follows:
\[
\rho(a) =
\begin{bmatrix}
\alpha & \alpha \beta \gamma + \alpha^{-1} \\
0 & \alpha^{-1}
\end{bmatrix},~~
\rho(b) =
\begin{bmatrix}
\gamma & 0 \\
-\beta^{-1}-\gamma^{-1} & \gamma^{-1}
\end{bmatrix}
\]
where $a$ and $b$ are generators of $\pi_1(P)$ which are the homotopy classes of the boundary components $A$ and $B$ respectively, and $\alpha, \beta, \gamma : \mathbb{R}_{>0}^3 \rightarrow \mathbb{R}_{>0}$ are defined by
\[
\alpha(l_A, l_B, l_C) = e^{l_A/2},~~
\beta(l_A, l_B, l_C)  = e^{(l_C - l_A)/2},~~
\gamma(l_A, l_B, l_C) =  e^{-l_B/2}.
\]
\end{proposition}

To use the Bonahon-Dreyer parametrization, we take the geodesic maximal lamination $\mathcal{L}$ on $P$ which consists of three biinfinite leaves $h_{AB}, h_{BC}$ and $h_{CA}$, and three closed leaves $A, B$ and $C$ (see Figure 5 in Section 5.1). 
The lamination $\mathcal{L}$ induces an ideal triangulation which separates $P$ into two ideal triangles $T_0$ and $T_1$. 
When we apply the Bonahon-Dreyer parameterization to our case, $\Phi_{\mathcal{L}}$ is the map defined by
{\small
\[ \Phi_{\mathcal{L}}([\rho]) = (\sigma_p^{\rho}(h_{AB}), \cdots, \sigma_p^{\rho}(h_{BC}), \cdots, \sigma_p^{\rho}(h_{CA}), \cdots, \tau_{pqr}^{\rho}(\widetilde{T}_0,v_0), \cdots, \tau_{pqr}^{\rho}(\widetilde{T}_1, v_1), \cdots) \] 
}
The following theorem is the main result of this paper which is an explicit description of ${\rm Fuch}_n(P)$ in the Bonahon-Dreyer's parameter space.
\begin{theorem}[Theorem 5.3]
The Bonahon-Dreyer coordinate $\Phi_{\mathcal{L}}([\rho_n])$ of any representation $[\rho_n] \in {\rm Fuch}_n(P)$ can be explicitly computed by using $\alpha, \beta$ and $\gamma$ which are the parameters of Fuchsian representations in Proposition1.1.  
\end{theorem} 
Especially we give general formulae of $\sigma_p^{\rho}(h_{AB}), \sigma_p^{\rho}(h_{BC}), \sigma_p^{\rho}(h_{CA}), \tau_{pqr}^{\rho}(\widetilde{T}_0,v_0)$ and $\tau_{pqr}^{\rho}(\widetilde{T}_1, v_1)$ to compute $\Phi_{\mathcal{L}}([\rho_n])$ for given $n$. The details of the computational results are in Theorem 5.3.

\subsection*{Structure of the paper}

This paper is organized as follows. 
In Section 2 we define Teichm\"uller spaces of surfaces, character varieties, and Hitchin components.
In particular we define Hitchin components of compact surfaces with boundary following \cite{BD2} and \cite{LM} in Section 2.3.
In Section 3, we prepare some invariants of flags and define Anosov representations.
Moreover we discuss the Anosov property of Hitchin representations which are important concepts to construct the Bonahon-Dreyer parameterization. 
In Section 4, we construct the Bonahon-Dreyer's parametrizations of Hitchin components.
In Section 4.4, we review Labourie-McShane's work \cite{LM} where they studied Hitchin representations of surfaces with boundary and defined Hitchin doubles which are extensions of Hitchin representations of surfaces with boundary to representations of doubled surfaces.
The Bonahon-Dreyer parametrizations of Hitchin components of surfaces with boundary are defined by using Hitchin doubles.
Finally, in Section 5, we compute Bonahon-Dreyer's parameters on Fuchsian loci. 
In Section 5.1, we parameterize Fuchsian representations of a pair of pants by hyperbolic lengths of the boundary components up to conjugacy.
Later, in Section 5.2, we compute the Bonahon-Dreyer parameters of Hitchin representations in Fuchsian loci.
We show that we can explicitly describe flag curves of representations in Fuchsian loci even in the cases of surfaces with boundary, and compute shearing and triangle invariants by using the parameters of Fuchsian representations.

\subsection*{Acknowledgement}
I would like to thank my supervisor, Ken'ichi Ohshika for his
warm encouragement. I also thank Hideki Miyachi who read my paper and gave me insightful comments and suggestions. 



\section{Teichm\"uller Spaces  and Hitchin components}
\subsection{Teichm\"uller Spaces and Hyperbolic structures}
Let $S$ be a compact, connected, oriented surface with negative Euler characteristic number $\chi (S)$. 
We denote by ${\rm Hyp}(S)$ the set of finite-volumed complete Riemannian metrics on $S$ with constant curvature $-1$ which make $\partial S$ totally geodesic if $\partial S$ is not empty, and denote by ${\rm Diff}_0(S)$ the identity component of the group of diffeomorphisms of $S$. 
The {\it Teichm\"uller space} of $S$ is the quotient space $\mathscr{T}(S) = {\rm Hyp}(S)/{\rm Diff}_0(S)$.
We call elements of $\mathscr{T}(S)$ {\it hyperbolic structures} of $S$.
Here we describe the Teichm\"uller space of a pair of pants. 
Insights of hyperbolic structures of geodesic hexagons imply that hyperbolic structures of a pair of pants are uniquely determined by the hyperbolic lengths of the boundary components. 
Let $P$ be a pair of pants whose boundary components are labeled as $A, B$ and $C$.
(See Figure 1.) 
We denote the length function associated to a hyperbolic structure $\textbf{m} \in \mathscr{T}(S)$ by $l_{\textbf{m}} : \mathcal{S} \rightarrow \mathbb{R}_{>0}$, where $\mathcal{S}$ is the set of simple closed curves on $S$. 
\begin{proposition}[\cite{IT}] 
The following map is a diffeomorphism.
\begin{equation*}
\mathscr{T}(P) \rightarrow \mathbb{R}_{>0}^3 ~:~ {\bf m} \mapsto (l_{{\bf m}}(A), l_{{\bf m}}(B), l_{{\bf m}}(C)).
\end{equation*}
\end{proposition}  

\begin{figure}[htbp]
\begin{center}
\includegraphics[width = 3cm, height=3cm]{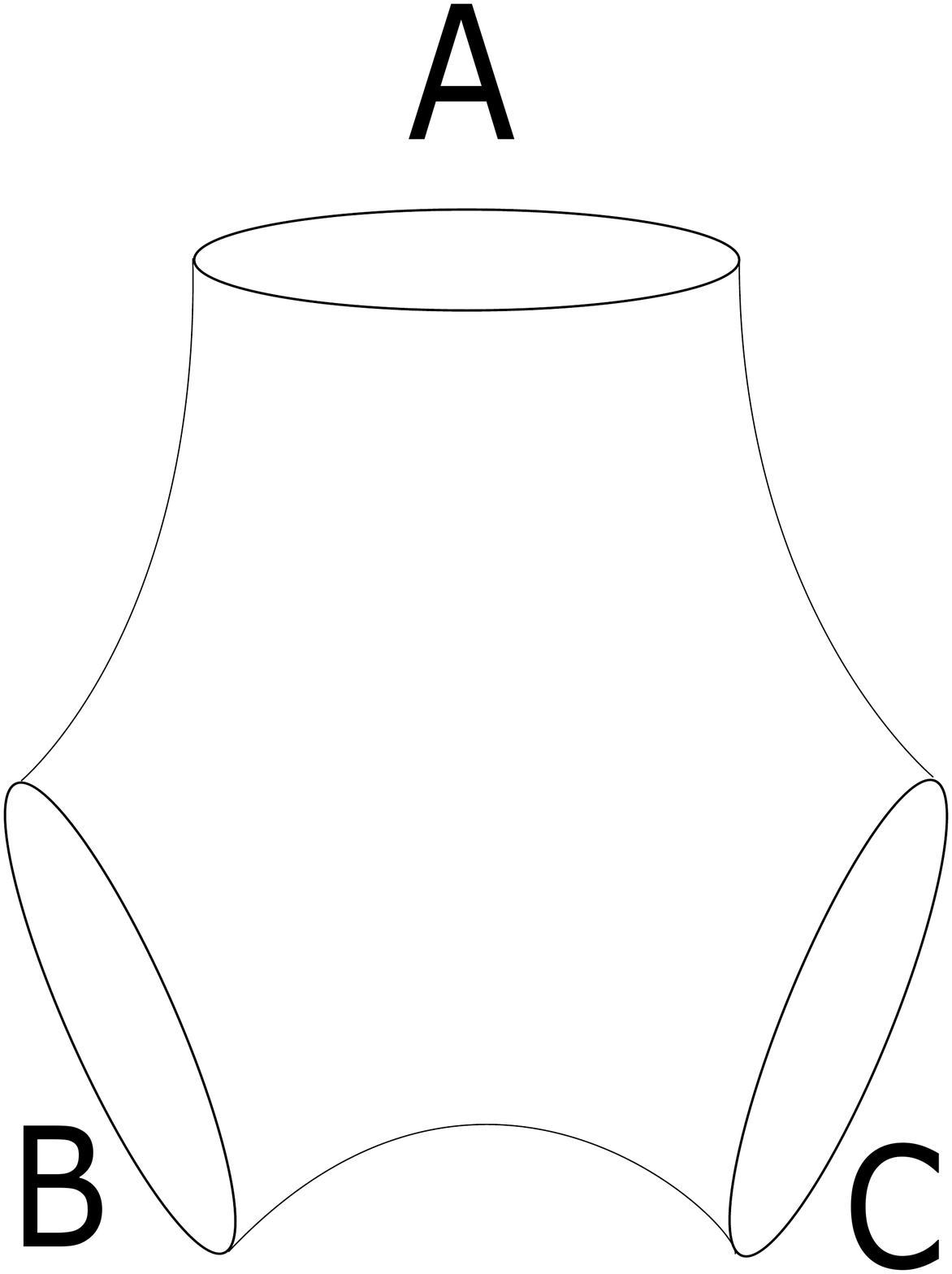}
\caption{A pair of pants.}
\label{figure 1}
\end{center}
\end{figure}

We can also consider hyperbolic structures from the viewpoint of {\it geometric structures}.
Let $X = \mathbb{H}^n$ be the upper half space model of the $n$-dimensional hyperbolic space, and $G  ={\rm Isom}^+(X)$ the group of orientation-preserving isomerties of $X$.   

\begin{definition}
Let $M$ be a manifold. 
An $(X,G)$-atlas on $M$ is an atlas $\{(U_i, \phi_i) \}_{i \in I}$ which is defined by the following conditions:
\begin{itemize}
\item[(1)] Each $U_i$ is an open subset of $M$ such that $M = \bigcup_iU_i$, and $\phi_i$ is an homeomorphism $\phi_i : U_i \rightarrow X$ which is onto the image.
\item[(2)] Any transition map $\phi_i \circ \phi_j^{-1} : \phi_j(U_i \cap U_j) \rightarrow\phi_i(U_i \cap U_j)$ is given by the restriction of an element in $G$. 
\end{itemize}
We say that a manifold $M$ has an $(X,G)$-structure if $M$ has a maximal $(X,G)$-atlas. 
A manifold with an $(X,G)$-structure is called an $(X,G)$-manifold. 
\end{definition}
Note that an $(X,G)$-manifold $M$ has a metric induced from the metric of $X$ through the $(X,G)$-atlas of $M$.
An $(X,G)$-manifold $M$ is said to be {\it complete} if $M$ is metrically complete with respect to the induced metric. 
We can check that $M$ is a complete $(X,G)$-manifold if and only if $M$ has a hyperbolic structure in the sense of Riemannian metrics.
Therefore we identify complete $(X,G)$-structures with hyperbolic structures, and we use the terms, {\it hyperbolic structures} and {\it hyperbolic manifolds} instead of complete $(X,G)$-structures and complete $(X,G)$-manifolds. 

Definition 2.2 gives us tools, called developing maps and holonomy representations, which are useful to study  geometric structures on manifolds. 
Let $M$ be a hyperbolic manifold with an $(X,G)$-atlas $\mathcal{U}=\{(U_i, \phi_i) \}_{i \in I}$, and $p_0 \in M$ a base point. 
A developing map is a map from the universal covering of $M$ to the model space $X$.
To construct this map, recall that the universal covering of $M$ can be constructed as the quotient space of pairs of paths and their terminal points;
\begin{equation*}
\widetilde{M} = \{ (C, p) ~|~ \text{$C$ is a path from $p_0$ and $p$ is the terminal point of $C$} \}/ \sim.
\end{equation*}  
where $(C_1,p_1) \sim (C_2, p_2)$ if $p_1 = p_2$, and $C_1$ is homotopic to $C_2$. 
For an equivalent class $[C,p] \in \widetilde{M}$, we fix a representative $(C,p)$.
Cover the curve $C$ by using a finite number of charts $(U_{i_1}, \phi_{i_1}), \cdots, (U_{i_n}, \phi_{i_n}) \in \mathcal{U}$ where $p_0 \in U_{i_1}$, $p \in U_{i_n}$, and $U_{i_j}$ and $U_{i_{j+1}}$ are adjacent for any $j$.
The transition map from $U_{i_2}$ to $U_{i_1}$ is the restriction of an element $g \in G$, and we can glue these charts $U_{i_1}$ and $U_{i_2}$ by $g$.
Repeating this operation, we construct a copy $\widetilde{C}$ of $C$ in $X$.
Then we define a map ${\bf dev} : \widetilde{M} \rightarrow X$ so that ${\bf dev}([C,p])$ equals to the terminal point of  $\widetilde{C}$. 
It is known that ${\bf dev}$ is well-defined, locally homeomorphic and unique up to the left-composition of isometries.
The map ${\bf dev}$ is called the {\it developing map}. 

Consider the covering transformation of $\pi_1(M)$ on the universal covering of an hyperbolic manifold $M$ with a developing map $\textbf{d}$. 
Then there exists a representation $\rho_{{\bf d}} : \pi_1(M) \rightarrow G$ so that the developing map $\textbf{d}$ is $(\rho_{{\bf d}}, \pi_1(M))$-equivariant. 
This representation is often called the {\it holonomy representation} of the fundamental group of $M$ associated to $\textbf{d}$.
The following proposition is well-known.
\begin{proposition}
For any developing map $\textbf{d}$ of an hyperbolic manifold $M$ and the holonomy representation $\rho_{{\bf d}}$ associated to $\textbf{d}$, the representation $\rho_{{\bf d}}$ is discrete and faithful. 
The discrete group $\rho_{{\bf d}}(\pi_1(M))$ acts on $\textbf{d}(X)$ properly-discontinuously and freely. 
Moreover $M$ is isometric to $\textbf{d}(X)/\rho_{{\bf d}}(\pi_1(M))$.
\end{proposition}
When we replace a developing map $\textbf{d}$ with the composition $g \cdot \textbf{d}$ for any $g \in G$, the corresponding holonomy representation changes to the conjugation $g \cdot \rho_{{\bf d}} \cdot g^{-1}$. 
Thus the holonomy representation of a hyperbolic manifold $M$ is uniquely determined up to conjugacy.

Note that when $n=2$, the isometry group $G$ is isomorphic to ${\rm PSL}_2(\mathbb{R})$, and note that any surface $S$ with $\chi(S) < 0$ has an $(\mathbb{H}^2, {\rm PSL}_2(\mathbb{R}))$-structure. 
Then deformation theory of hyperbolic structures can be translated to representation theory of $\pi_1(S)$ into ${\rm PSL}_2(\mathbb{R})$. 
In particular Teichm\"uller spaces are identified with spaces of representations. We remark that one can define the following well-defined map which sends hyperbolic structures to associated holonomy representations.
\[  {\rm Hol} : \mathscr{T}(S) \rightarrow {\rm Hom}(\pi_1(S), {\rm PSL}_2(\mathbb{R}))/{\rm PSL}_2(\mathbb{R})-\mbox{conjugacy}. \]
We call this map the {\it holonomy map}. 
Spaces of representations are well studied in many contexts. 
In the next section we review a theory of representation varieties of surface groups and their topologies.

\subsection{Representation varieties for surface groups}

The ${\rm PSL}_n(\mathbb{R})$-representation variety of the surface $S$ is the space of representations $\mathcal{R}_n(S)={\rm Hom}(\pi_1(S), {\rm PSL}_n(\mathbb{R}))$.
The representation variety $\mathcal{R}_n(S)$ is regarded as an algebraic set.
In fact, if we fix a presentation of $\pi_1(S)$ with $k$-generators, then $\mathcal{R}_n(S)$ is identified with a subset in the product ${\rm PSL}_n(\mathbb{R})^k$, and the relations of $\pi_1(S)$ give defining polynomials of $\mathcal{R}_n(S)$.
The Lie group ${\rm PSL}_n(\mathbb{R})$ acts on the representation variety by conjugation.
The quotient $\mathcal{X}_n(S)=\mathcal{R}_n(S)// {\rm PSL}_n(\mathbb{R})$ by the action is called the ${\rm PSL}_n(\mathbb{R})$-{\it character variety} of $S$, where ``$//$'' means the GIT quotient which is a notion of algebraic geometry. (See \cite{S} for details.)

For the rest of this section, we assume that $S$ is closed and has the genus $g \geq 2$ unless stated otherwise. 
In \cite {Go}, Goldman detected the number of connected components of $\mathcal{R}_2(S)$.
We say that a discrete faithful representation $\rho$ is orientation-preserving if the natural equivalence relation $S \rightarrow \mathbb{H}^2/\rho(\pi_1(S))$ has degree $+1$.  
\begin{theorem}[Goldman\cite{Go}]
The number of connected components of $\mathcal{R}_2(S)$ equals to $4g-3$.
Moreover the subset ${\rm Fuch}_2(S)$ of $\mathcal{X}_2(S)$ consisting of orientation-preserving discrete faithful representations is a connected component of $\mathcal{X}_2(S)$.  
\end{theorem}
In particular, $\mathscr{T}(S)$ is diffeomorphic to ${\rm Fuch}_2(S)$ by the holonomy map ${\rm Hol}: \mathscr{T}(S) \rightarrow \mathcal{X}_2(S)$.
In this reason the component ${\rm Fuch}_2(S)$ is called the {\it Teichm\"uller component}. 
A representation $\rho \in \mathcal{R}_2(S)$ is called a {\it Fuchsian representation} if $[\rho] \in {\rm Fuch}_2(S)$. 
We denote the subset of Fuchsian representations in $\mathcal{R}_2(S)$ by $\mathcal{F}_2(S)$.

\subsection{Hitchin components}
We define Hitchin components of closed surfaces. 
First we recall Hitchin's work about the number of components of $\mathcal{X}_n(S)$.
Hitchin components were originally defined by Hitchin who developed the above result of Goldman by using Higgs bundle theory.
\begin{theorem}[Hitchin\cite{Hi}]
For $n \geq 3$, $\mathcal{X}_n(S)$ has three connected components if $n$ is odd, and has six connected components if $n$ is even.
\end{theorem}
Consider the irreducible representation $\iota_n: {\rm PSL}_2(\mathbb{R}) \rightarrow {\rm PSL}_n(\mathbb{R})$ which is given by the irreducible representation of ${\rm SL}_2(\mathbb{R})$ on the symmetric power $Sym^{n-1}(\mathbb{R}^2)$. 
(See \cite{Hi}.)
This representation induces the map $(\iota_n)_* : \mathcal{R}_2(S) \rightarrow \mathcal{R}_n(S)$ which is defined by the left composition, and induces the well-defined map $(\iota_n)_* : \mathcal{X}_2(S) \rightarrow \mathcal{X}_n(S)$.
\begin{definition}
The Hitchin component of $S$ is the connected component of $\mathcal{X}_n(S)$ which contains $(\iota_n)_*({\rm Fuch}_2(S))$. We denote it by ${\rm Hit}_n(S)$.  
\end{definition}

Following a formulation of Labourie-McShane \cite{LM}, we define Hitchin components of surfaces with boundary. 
Let $S_{\partial}$ be a surface with no punctures, boundaries $C_1, \cdots, C_b~(b>0)$ and $\chi(S_{\partial}) <0$. 
In this case we must restrict images of boundary curves. 
An element $A \in {\rm PSL}_n(\mathbb{R})$ is called {\it purely-loxodromic} if $A$ is conjugate to an element in the interior of the Weyl chamber, that is, $A$ is a real split with distinct eigenvalues. 
We call $\rho \in \mathcal{R}_n(S_{\partial})$ {\it boundary-loxodromic} if for $[C_i] \in \pi_1(S_{\partial})~(i=1, \cdots, b)$, $\rho([C_i])$ is a purely-loxodromic element. 
We set the subset $\mathcal{R}_n^{loxo}(S_{\partial})$ of $\mathcal{R}_n(S_{\partial})$ by
\[  \mathcal{R}_n^{loxo}(S_{\partial}) = \{ \rho \in \mathcal{R}_n(S_{\partial}) ~|~ \rho \mbox{ is boundary-loxodromic}  \}. \] 
Then we define the Hitchin component ${\rm Hit}_n(S_{\partial})$ by the component of the quotient $\mathcal{X}_n^{loxo}(S_{\partial})= \mathcal{R}_n^{loxo}(S_{\partial}) //  {\rm PSL}_n(\mathbb{R})$ which contains $(\iota_n)_*({\rm Fuch}_2(S_{\partial}))$.
Here note that metrics in Teichm\"uller spaces make the boundary totally-geodesic, the Teichm\"uller components ${\rm Fuch}_2(S_{\partial})$ are contained in $\mathcal{X}_2^{loxo}(S_{\partial})$.

We prepare several terms. 
A representation $\rho \in \mathcal{R}_n(S)$ is called a {\it Hitchin representation} if $[\rho] \in {\rm Hit}_n(S)$. 
We denote the subset of Hitchin representations by $\mathcal{H}_n(S)$. 
The subset ${\rm Fuch}_n(S) = (\iota_n)_*({\rm Fuch}_2(S))$ is called the {\it Fuchsian locus} and a representation $\rho_n \in \mathcal{R}_n(S)$ with $[\rho_n] \in {\rm Fuch}_n(S)$ is called an {\it $n$-Fuchsian representation}. 
Any $n$-Fuchsian representation is of the form $\iota_n \circ \rho$ for some $\rho \in \mathcal{F}_2(S)$. 
The subset $\mathcal{F}_n(S)$ in $\mathcal{R}_n(S)$ denotes the set of $n$-Fuchsian representations.


\section{The Anosov property of Hitchin representations}

In this section we introduce the concept of Anosov representations following \cite{GW1}, and their application to Hitchin components. 
Here we consider only the case of closed surfaces.
The non-closed case is going to be discussed in Section 4.4. 

\subsection{Anosov representations}

Let $\Gamma= \pi_1(S)$ and $G$ a semisimple Lie group. 
Fix a hyperbolic structure of $S$. 
We consider the geodesic flow $\phi_t$ on the unit tangent bundle $T^1S$, and lift $\phi_t$ to the flow $\widetilde{\phi_t}$ on the unit tangent bundle $T^1\tilde{S}$ over the universal covering of S associated to the hyperbolic structure of $S$. 
Let $(P^+, P^-)$ be an opposite pair of parabolic subgroups of $G$ where $P^+$ and $P^-$ are called {\it opposite} if their intersection $L=P^+ \cap P^-$ is reductive. 
The quotient $G/L$ is contained in the product $G/P^+ \times G/P^-$ as the unique orbit space, and its tangent space $T_x G/L$ at a point $x \in G/L$ can be decomposed into two directions $E_x^+ \oplus E_x^-$ along the decomposition of the tangent space of the product space $T_x(G/P^+ \times G/P^-)=T_{x^+}G/P^+ \oplus T_{x^-}G/P^-$ where $(x^+, x^-)$ is a point of $G/P^+ \times G/P^-$ corresponding to the point $x$. 
Let $E^+$ and $E^-$ be the distribution consisting of $E_{x^+}^+$ and $E_{x-}^-$ over $G/L$ respectively. 
Now we consider the associated bundle $T^1\tilde{S} \times_{\rho} G/L$ over the unit tangent bundle $T^1S$ for a representation $\rho: \Gamma \rightarrow G$.
\begin{definition}
A representation $\rho : \Gamma \rightarrow G$ is said to be $(P^+, P^-)$-Anosov if 
\begin{itemize}
\item[(1)] There exists a section $\sigma : T^1S \rightarrow T^1\tilde{S} \times_{\rho} G/L$ which is locally constant along the flow $\phi_t$.
\item[(2)]The lifted action of $\phi_t$ on the pull-back $\sigma^*E^+$ (resp. $\sigma^*E^-$) satisfies the dilating (resp. contracting) property.
\end{itemize}
\end{definition}
In Definition 3.1, the condition (2) means that there exists a continuous family of norms $\{ ||\cdot||_p \}_{p \in T^1S}$ on the fibers of the distribution $\sigma^*E^+$ (resp. $\sigma^*E^-$) and positive constants $a, A >0$ such that  
\begin{equation*}
|| \phi_{-t}(v) ||_{\phi_{-t}(p)} \leq A e^{-at} ||v||_{p},~~ (\mbox{resp. }|| \phi_{t}(v) ||_{\phi_{t}(p)} \leq A e^{-at} ||v||_{p})
\end{equation*}
for any $t \geq 0$ and any $v \in \sigma^*E^+(\mbox{resp. } \sigma^*E^-)$ which is in the fiber on an arbitrary point $p \in T^1S$. 
The section $\sigma$ in Definition 3.1 is called the {\it Anosov section}.
The Anosov section induces a map from $\partial_{\infty}\tilde{S}$ to $G/P^+$ and $G/P^-$ as follows. 
Note that the orbits of the geodesic flow on $T^1\tilde{S}$ are identified with doubles in $\partial_{\infty}\tilde{S} \times \partial_{\infty}\tilde{S}$, and that the set of oriented geodesic leaves on $\tilde{S}$ is identified with the set $(\partial_{\infty}\tilde{S})^{(2)} = \partial_{\infty}\tilde{S} \times \partial_{\infty}\tilde{S} \setminus \Delta$ where $\Delta$ is the diagonal set. 
Since the Anosov section $\sigma$ is locally constant, the pull-back $\widehat{\sigma} : T^1\tilde{S} \rightarrow G/L$ is $\widetilde{\phi}_t$-invariant.
Hence $\widehat{\sigma}$ induces a map defined on the quotient with the flow action $T^1\tilde{S}/(\widetilde{\phi}_t)= (\partial_{\infty}\tilde{S})^{(2)}$. 
We denote the map by
\begin{equation*}
\widehat{\sigma}=(\xi^+, \xi^-) : (\partial_{\infty}\tilde{S})^{(2)} \rightarrow G/L \subset G/P^+ \times G/P^-.
\end{equation*}
The map $\widehat{\sigma}$ factors through the projection from $ (\partial_{\infty}\mathbb{H}^2)^{(2)}$ onto each factor, and gives us maps $\xi^{\pm} : \partial_{\infty}\tilde{S} \rightarrow G/P^{\pm}$. 
It is known that these maps are unique, and $\xi^+ = \xi^-$ under the identification of $G/P^+$ and $G/P^-$, so we denote by $\xi = \xi^+ = \xi^-$. 
(See \cite{GW1}.) 
The map $\xi$ is called the {\it Anosov map}, the {\it flag curve} or the {\it limit curve}. 

It is known that Hitchin representations of surface groups into $G={\rm PSL}_n(\mathbb{R})$ are $(B^+, B^-)$-Anosov where $B^+$ and $B^-$ are opposite Borel subgroups of $G$. (See \cite{L})
In this case $G/B^{\pm}$ are the flag manifolds ${\rm Flag}(\mathbb{R}^n)$.
Labourie and Fock-Goncharov studied properties of flag curves in the case of ${\rm PSL}_n(\mathbb{R})$, so-called the hyperconvexity and the positivity. 
Their results are summarized in Section 3.3 after defining some notions with flags.

\subsection{Projective invariants}
We define projective invariants of tuples of flags which play an important role in the Bonahon-Dreyer parameterization. 

A {\it flag} in $\mathbb{R}^n$ is a sequence of nested vector subspaces of $\mathbb{R}^n$
\[ F = ( \{0\} = F^{(0)} \subset F^{(1)} \subset F^{(2)} \subset \cdots \subset F^{(n)} = \mathbb{R}^n )\]
where ${\rm dim}F^{(i)} = i$. 
The set of these flags is called the {\it flag manifold}.
We denote it by ${\rm Flag}(\mathbb{R}^n)$.
A {\it generic} tuple of flags is a tuple $(F_1, F_2, \cdots, F_k)$ of a finite number of flags $F_1, F_2, \cdots, F_k \in {\rm Flag}(\mathbb{R}^n)$ such that if $n_1, \cdots, n_k$ are nonnegative integers satisfying $n_1 + \cdots + n_k = n$, then $F_1^{(1)} \cap \cdots \cap F_k^{(n_k)} = \{ 0 \}$. 

Now we define invariants called {\it triple ratios} and {\it double ratios} of
 generic tuples of flags. 
The triple ratio is defined for generic triples of flags as follows. 
Let $(E, F, G)$ be a generic triple of flags, and $p,q,r \geq 1$ integers satisfying $p+q+r = n$. 
Fix nonzero elements $e^{(i)}, f^{(i)}$ and $g^{(i)}$ which are in the wedge products $\bigwedge^{(i)}E^{(i)}, \bigwedge^{(i)}F^{(i)}$ and $\bigwedge^{(i)}G^{(i)}$ for $i=1, \cdots , n$ respectively. 
We denote the wedge product $e^{(p)} \wedge f^{(q)} \wedge g^{(r)}$ of chosen elements $e^{(p)}, f^{(q)}$ and $g^{(r)}$ by $X(p,q,r)$, which is an element in $\bigwedge^{(n)}\mathbb{R}^n$ since $p+q+r=n$.
\begin{definition}
The $(p,q, r)$-th triple ratio $T_{pqr}(E,F,G)$ for a generic triple of flags $(E, F, G)$ is defined by
\[ T_{pqr}(E,F,G) = \dfrac{X(p+1, q, r-1)}{X(p-1, q, r+1)} \cdot \dfrac{X(p, q-1, r+1)}{X(p, q+1, r-1)} \cdot \dfrac{X(p-1, q+1, r)}{X(p+1, q-1, r)}. \]
\end{definition}
Triple ratios are valued in $\mathbb{R}$ under an identification $\bigwedge^{(n)}\mathbb{R}^n \cong \mathbb{R}$. 
This is well-defined for a choice of elements $e^{(i)}, f^{(i)}$ and $g^{(i)}$, and an identification $\bigwedge^{(n)}\mathbb{R}^n \cong \mathbb{R}$. 
If the index of $X(p,q,r)$ contains $0$ then we ignore the corresponding terms, i.e. for example, $e^{(0)} \wedge f^{(q)} \wedge g^{(n-q)} =  f^{(q)} \wedge g^{(n-q)}$. 
For replacements of an ordering of a generic tuple $(E,F,G)$, triple ratios well behave.
\begin{proposition}[\cite{BD2}]
For a generic tuple of flags $(E,F,G)$, 
\[ T_{pqr}(E,F,G) = T_{qrp}(F,G,E) = T_{qpr}(F,E,G)^{-1}. \]
\end{proposition}

Next we define double ratios. 
Let $(E,F,G,G')$ be a generic quadruple of flags, and $p$ an integer with $1 \leq p \leq n-1$.
We choose nonzero elements $e^{(i)}, f^{(i)}, g^{(i)}$ and $g'^{(i)}$ in $\bigwedge^{(i)}E^{(i)}, \bigwedge^{(i)}F^{(i)}, \bigwedge^{(i)}G^{(i)}$ and $\bigwedge^{(i)}G'^{(i)}$ respectively.
Let us set $Y(i) = e^{(i)} \wedge f^{(n-i-1)} \wedge g^{(1)}$ and $Y'(i) = e^{(i)} \wedge f^{(n-i-1)} \wedge g'^{(1)}$.
\begin{definition}
The $p$-th double ratio $D_p(E,F,G,G')$ is defined by
\[ D_p(E,F,G,G') = - \dfrac{Y(p)}{Y'(p)} \cdot \dfrac{Y'(p-1)}{Y(p-1)}. \]
\end{definition}
This is also well-defined and valued in $\mathbb{R}$ as triple ratios. 
We remark that triple ratios and double ratios are invariant under the action of projective automorphisms.

\subsection{The Anosov property}
We make in the summary results of Labourie and Fock-Goncharov. 
Recall that  any Hitchin representation $\rho \in \mathcal{H}_n(S)$ is Anosov, so there exists the corresponding flag curve $\xi_{\rho}$. 
The flag curves of Hitchin representations depend on the eigenspaces of matrices in $\rho(\pi_1(S))$. 
Now we recall a property of eigenvalues of Hitchin representations.
\begin{proposition}[Labourie \cite{L}, Bonahon-Dreyer \cite{BD2}]
Let $\rho \in \mathcal{H}_n(S)$, and $\gamma$ be a nontrivial element of $\pi_1(S)$. Then $\rho(\gamma)$ has a lift $ \widetilde{\rho(\gamma)}  \in {\rm SL}_n(\mathbb{R})$ whose eigenvalues are distinct and positive.
\end{proposition}
For any nontrivial element $\gamma \in \pi_1(S)$, we denote the eigenvalues of a lift $\widetilde{\rho(\gamma)}$ by
\[  \lambda^{\rho}_1(\gamma) > \lambda^{\rho}_2(\gamma) >  \cdots > \lambda^{\rho}_n(\gamma) >0. \]
Let $L_k$ be the one dimensional  eigenspace of $\lambda_k^{\rho}(\gamma)$. 
We associate to $\rho(\gamma)$ the two flags $E$ and $F$ in ${\rm Flag(\mathbb{R}^n)}$ defined by
\[ E^{(i)} = \oplus_{k=1}^i L_k \mbox{ and } F^{(i)} = \oplus_{k=n-i+1}^n L_k \]
which are called the {\it stable flag} and the {\it unstable flag} of $\rho(\gamma)$ respectively.

Fix a hyperbolic structure on $S$ and the associated universal covering $\widetilde{S}$ with the visual boundary $\partial_{\infty}\widetilde{S}$. 
\begin{theorem}[Labourie \cite{L}, Fock-Goncharov \cite{FG}]
For any $\rho \in \mathcal{H}_n(S)$, there exists a unique continuous map $\xi_{\rho} : \partial_{\infty}\widetilde{S} \rightarrow {\rm Flag}(\mathbb{R}^n) $ satisfying the following conditions.
\begin{itemize}
\item[(1)] For the attracting point $\gamma^+$ of a nontrivial $\gamma \in \pi_1(S)$, the flag $\xi_{\rho}(\gamma^+)$ is the stable flag of $\rho(\gamma)$.
\item[(2)]The map $\xi_{\rho}$ is $\rho$-equaivariant for $\pi_1(S)$-action on $\partial_{\infty}\widetilde{S}$, that is, for any element $\gamma \in \pi_1(S)$, it holds that $\rho(\gamma) \cdot \xi_{\rho}(x) = \xi_{\rho}(\gamma \cdot x)$.
\item[(3)] For any $x_1 \cdots, x_k \in \partial_{\infty}\widetilde{S}$, a tuple $(\xi_{\rho}(x_1), \cdots,  \xi_{\rho}(x_k))$ is generic.
\item[(4)] For any triple $(x, y, z)$ of distinct points $x,y,z \in \partial_{\infty}\widetilde{S}$, triple ratios of a generic triple $(\xi_{\rho}(x), \xi_{\rho}(y), \xi_{\rho}(z))$ are positive.
\item[(5)] For any quadruple $(x, z, y, z')$ of distinct points $x,z,y,z' \in \partial_{\infty}\widetilde{S}$ which are in counterclockwise ordering, double ratios of $(\xi_{\rho}(x), \xi_{\rho}(y), \xi_{\rho}(z), \xi_{\rho}(z'))$ are positive.
\end{itemize}
\end{theorem}


\section{The Bonahon-Dreyer parameterization}
We review the Bonahon-Dreyer parameterization in this section. 
Through this section, we suppose that $S$ is closed and has the genus $g \geq 2$ unless otherwise indicated..
Moreover we fix a hyperbolic structure of $S$. 
Details of this section are in \cite{BD2}. 
\subsection{Short arcs for laminations}
To define their parameterization, we give more additional data to maximal geodesic laminations.
Let $\mathcal{L}$ be a maximal geodesic lamination on $S$.
We suppose the orientation of each leaf of $\mathcal{L}$, where orientations of leaves are given independently.  
Moreover we fix certain short transverse arcs for all closed leaves.
For a closed leaf $g$ of $\mathcal{L}$, we fix a transverse arc $k$ so that the arc $k$ satisfies the following conditions:
\begin{itemize}
\item[(1)] The arc $k$ is transverse to $\mathcal{L}$ and the intersection $k \cap g$ is just one point $x$. 
\item[(2)] Let $k'$ and $k''$ be arcs given by cutting $k$ by $g$. Then, for each component $c = k'$ or $k''$, there exists an immersion $f : c \times [0, +\infty) \rightarrow S$ such that (i)$f(x,0) = x$ and (ii) for any $x \in \mathcal{L} \cap c$, the subset $\{x \} \times [0, \infty)$ parametrizes with unit speed  a geodesic spiraling along $g$ and the image $f(x,[0, \infty))$ is contained in some leaf of $\mathcal{L}$.
\end{itemize}

\subsection{Invariants of Hitchin representations}
We define two kinds of invariants of Hitchin representations, called the {\it triangle invariant} and the {\it shearing invariant}. 
We denote by $\widetilde{S}$ the universal covering of $S$ associated to the fixed hyperbolic structure.
To construct the invariants, we fix a maximal geodesic lamination $\mathcal{L}$ on $S$ which consists of closed leaves $g_1, \cdots, g_s$ and biinfinite leaves $h_1, \cdots, h_t$ which have independent orientations. 
Note that $\mathcal{L}$ induces an ideal triangulation $\Delta$ of $S$. 
We assume an existence of a short transverse arc $k_i$ for each closed leaf $g_i$, which satisfies the conditions in the previous subsection. 
Let $\rho$ be a Hitchin representation in $\mathcal{H}_n(S)$ and $\xi_{\rho}$ the associated flag curve.

The triangle invariant is defined for the ideal triangles $T_1, \cdots, T_u$ given by the ideal triangulation $\Delta$. 
Fix a lift $\widetilde{T_i}$ in $\widetilde{S}$ of $T_i$ and an ideal vertex $v_i$ of $\widetilde{T_i}$ arbitrarily. 
We label the other ideal vertices as the symbols $v'_i, v''_i$ so that $v_i, v'_i, v''_i$ are in clockwise ordering.
(See Figure 2.) 
Let $p,q,r$ be positive integers with $p+q+r=n$.
\begin{definition}
The $(p,q,r)$-th triangle invariant $\tau^{\rho}_{pqr}(\widetilde{T_i}, v_i)$ of a Hitchin representation $\rho$ for $(\widetilde{T_i}, v_i)$ is defined by
\[ \tau_{pqr}^{\rho}(\widetilde{T_i}, v_i) = \log T_{pqr}(\xi_{\rho}(v_i), \xi_{\rho}(v'_i), \xi_{\rho}(v''_i)). \]
\end{definition}

Triangle invariants are independent of a choice of the lift $\widetilde{T_i}$ since flag curves of Hitchin representations are $\rho$-equivariant, so  we denote the triangle invariant by $\tau^{\rho}_{pqr}(T_i, v_i)$ simply.
For replacements of a choice of vertex $v_i$, triangle invariants satisfy the following property. 
\begin{proposition}[\cite{BD2}]
Fix a lift $\widetilde{T_i}$ of the ideal triangle $T_i$. Let $v_i, v'_i, v''_i$ be vertices of $\widetilde{T_i}$ which are in clockwise ordering. Then the following relation holds.
\[  \tau_{pqr}^{\rho}(T_i, v_i) = \tau_{qrp}^{\rho}(T_i, v'_i) = \tau_{rpq}^{\rho}(T_i, v''_i). \]
\end{proposition}

\begin{figure}[htbp]
\begin{minipage}{0.45\hsize}
\begin{center}
\includegraphics[width = 40mm, height=4cm]{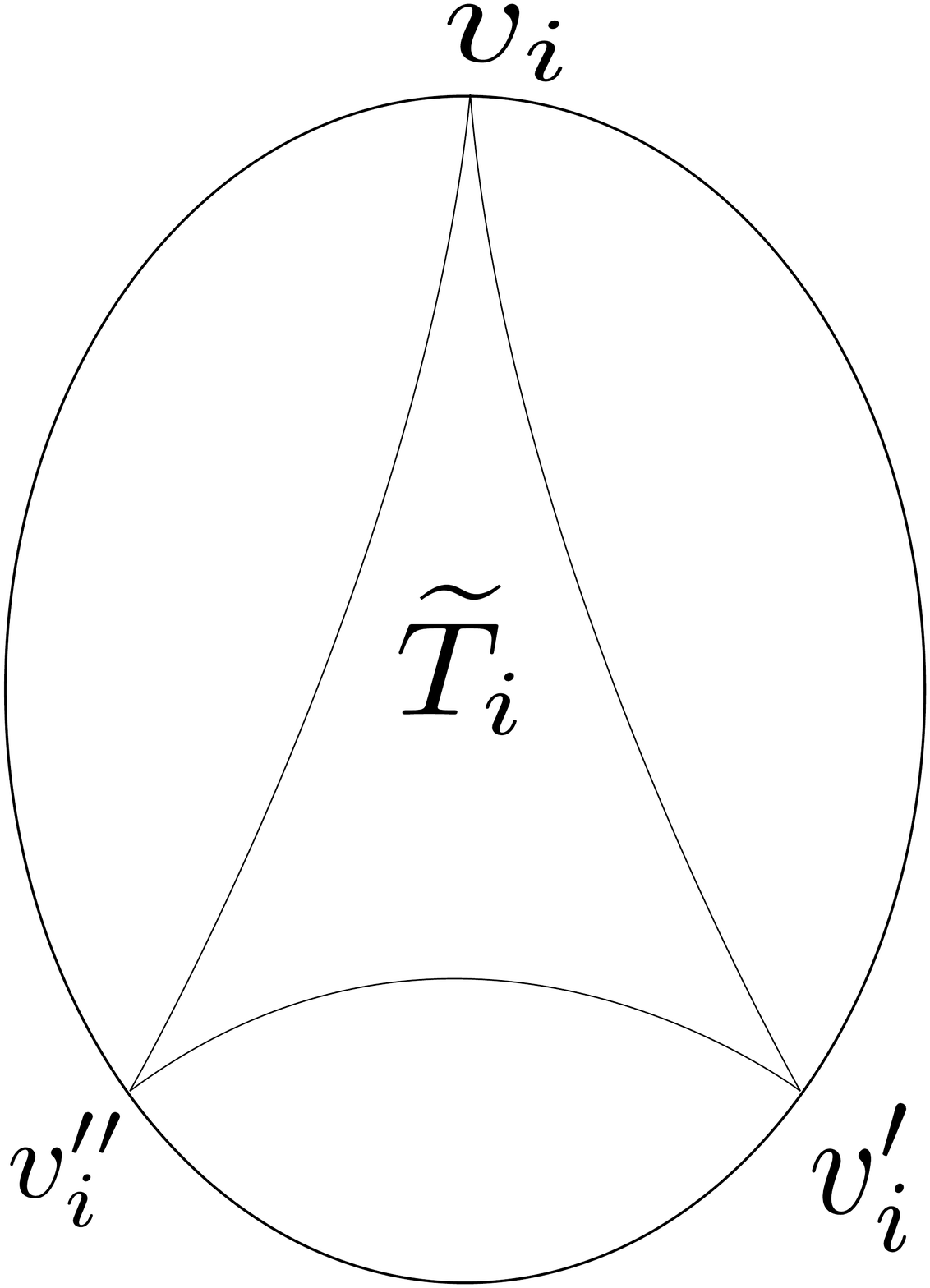}
\caption{}
\label{figure 2}
\end{center}
\end{minipage}
\begin{minipage}{0.45\hsize}
\begin{center}
\includegraphics[width=40mm, height=4cm]{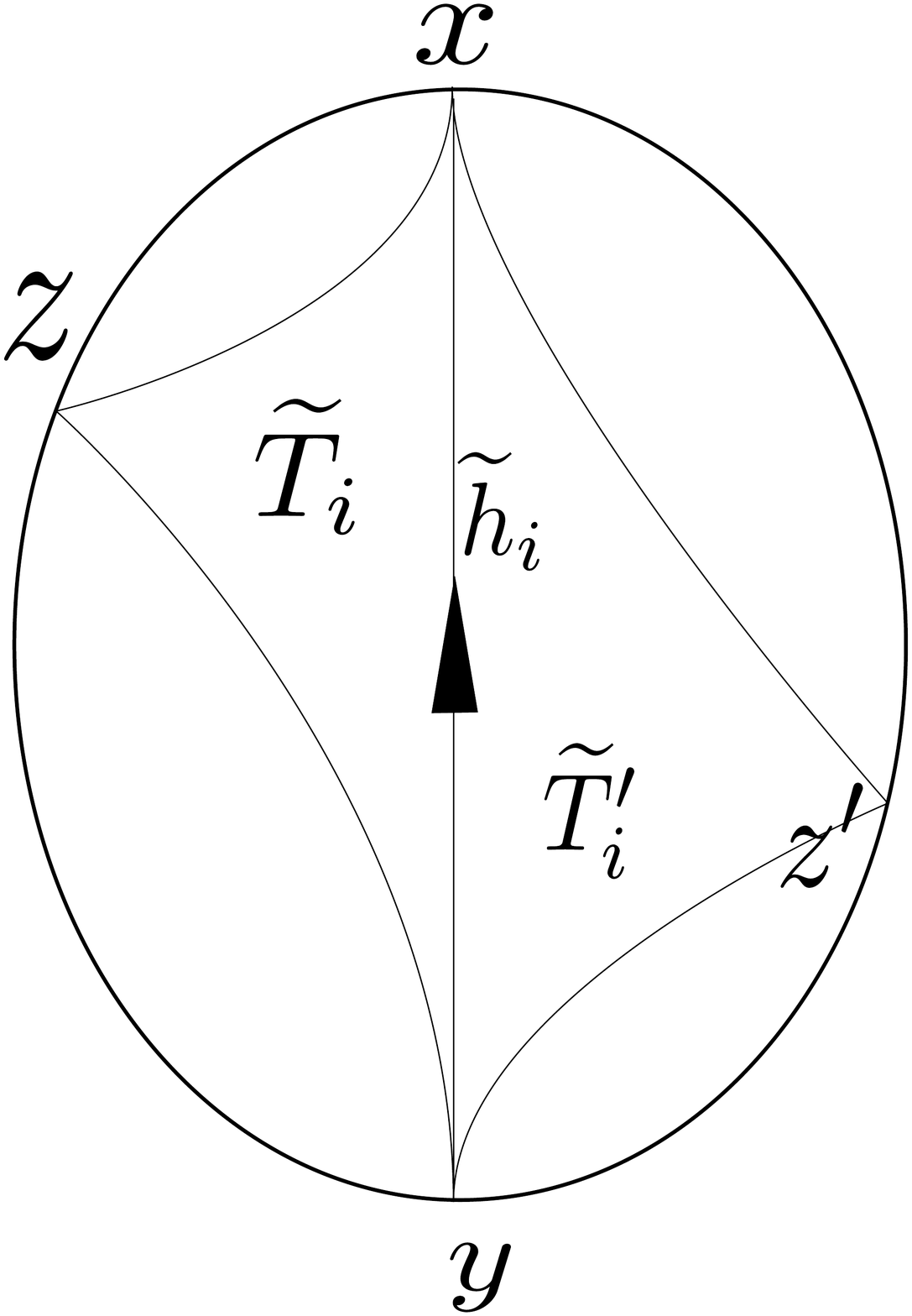}
\caption{}
\label{figure 3}
\end{center}
\end{minipage}
\end{figure}

The shearing invariant is defined for each leaf of $\mathcal{L}$.
At first we define the shearing invariant of a biinfinite leaf $h_i$. 
In the ideal triangulation $\Delta$, the leaf $h_i$ is adjacent to two ideal triangles. 
Let $T_i$ and $T'_i$ be such triangles which are on the left and right side with respect to the orientation of $h_i$. 
We lift $h_i$ to a geodesic $\widetilde{h_i}$ in $\widetilde{S}$, and also lift $T_i$ and $T_i'$ to two ideal triangles $\widetilde{T_i}$ and $\widetilde{T'_i}$ such that $\widetilde{h_i}$ is an edge of $\widetilde{T_i}$ and $\widetilde{T'_i}$ as Figure 3. 
We denote the starting and terminal point of $\widetilde{h_i}$ by $y$ and $x \in \partial_{\infty}\widetilde{S}$, and the other vertices of $\widetilde{T}_i$ and $\widetilde{T}'_i$ by $z$ and $z'$ respectively. 
Let $p$ be an integer with $1 \leq p \leq n-1$. 
Then we define the shearing invariant associated to $h_i$ which is independent of a choice of the lift $\widetilde{h_i}$. 
\begin{definition}
The $p$-th shearing invariant of a Hitchin representation $\rho$ for a biinfinite leaf $h_i$ is defined by 
\[ \sigma_p^{\rho}(h_i) = \log D_p(\xi_{\rho}(x), \xi_{\rho}(y), \xi_{\rho}(z), \xi_{\rho}(z')). \]
\end{definition}

In the case of closed leaves, we use the transverse arcs to take surrounding ideal triangles as in the case of biinfinite leaves. 
Fix a closed leaf $g_i$ with the transverse arc $k_i$ and consider the endpoints of $k_i$. 
Let $T_i$ and $T'_i$ be ideal triangles defined by $\Delta$ which are spiraling along $g_i$ such that $T_i$ (resp. $T'_i$) is on the left (resp. right) with respect to the orientation of $g_i$ and contains the endpoint of $k_i$. 
Lift  $g_i, k_i, T_i$ and $T_i'$ to a geodesics $\widetilde{g_i}$, an arc $\widetilde{k}_i$ and ideal triangles $\widetilde{T}_i$ and $\widetilde{T}'_i$ in $\widetilde{S}$ so that $\widetilde{k}_i$ intersects $\widetilde{g}_i$, and $\widetilde{T}_i$ and $\widetilde{T}'_i$ contain the endpoints of $\widetilde{k}_i$ as Figure 4. 
Label the starting point and the terminal point of the oriented geodesic $\widetilde{g}_i$ as $y$ and $x$ respectively. 
Now we can choose the edges $e_i$ and $e_i'$ of $\widetilde{T}_i$ and $\widetilde{T}'_i$ so that they intersect the endpoints of $\widetilde{g}_i$, and $\widetilde{g}_i, \widetilde{T}_i$ and $\widetilde{T}'_i$ are contained in the region which is bounded by $e_i$ and $e_i'$. 
Let  $z$ and $z'$ be the endpoints of $e_i$ and $e_i'$ other than $x$ and $y$, and $p$ an integer with $1 \leq p \leq n-1$.
\begin{definition}
The $p$-th shearing invariant of a Hitchin representation $\rho$ for a closed leaf $g_i$ is defined by 
\[ \sigma_{p}^{\rho}(g_i) = \log D_p(\xi_{\rho}(x), \xi_{\rho}(y), \xi_{\rho}(z), \xi_{\rho}(z')). \]
\end{definition}
Note that the shearing invariant $\sigma_{p}^{\rho}(g_i)$ are independent of a choice of the lift $\widetilde{g}_i$ and $\widetilde{k}_i$ since the flag curve $\xi_{\rho}$ is $\rho$-equivariant, and the image of the endpoints $\xi_{\rho}(x)$ and $\xi_{\rho}(y)$ are the stable and unstable flags of $\rho(\gamma)$ respectively. 

The triangle invariant and the shearing invariant are well-defined on Hitchin components, i.e. these invariants are independent of representatives of elements of Hitchin components since there is the relation $\xi_{g \cdot \rho \cdot g^{-1}} = g \cdot \xi_{\rho}$ for any $g \in {\rm PSL}_n(\mathbb{R})$.

\begin{figure}[htbp]
\begin{center}
\includegraphics[width = 4cm,height=4cm]{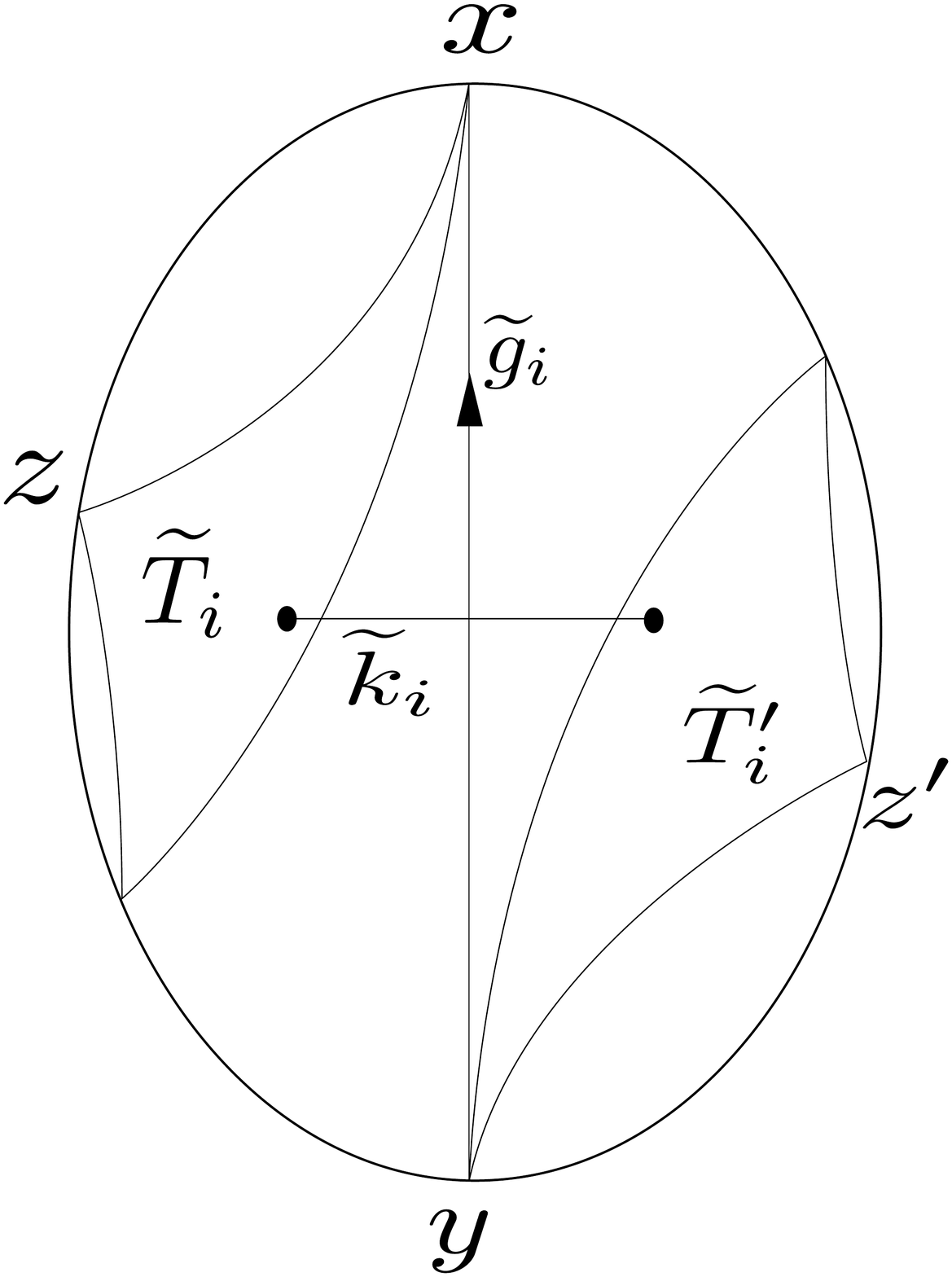}
\caption{}
\label{figure 4}
\end{center}
\end{figure}

\subsection{The construction of the parameter space}
The parameter space of ${\rm Hit}_n(S)$ given by Bonahon-Dreyer is the interior of a convex polytope defined by equations in a Euclidian space.
To explain the defining equations, we prepare length functions for Hitchin representations $\rho \in \mathcal{H}_n(S)$. 
Let $\lambda_1^{\rho}(\gamma) > \cdots > \lambda_n^{\rho}(\gamma)>0$ be the eigenvalues of the lift of $\rho(\gamma)$ in Proposition 3.5 for nontrivial $\gamma \in \pi_1(S)$. 
We define the {\it $k$-th length function} of the Hitchin representation $\rho$ by
\[ l_k^{\rho}(\gamma) = \log \dfrac{\lambda_k^{\rho}(\gamma)}{\lambda_{k+1}^{\rho}(\gamma)} \]
where $k=1, \cdots, n-1$. 
This function is independent of representatives of elements of ${\rm Hit}_n(S)$. 

There is a nice relation between length functions, triangle invariants and shearing invariants. 
Let $g$ be a closed leaf of the fixed lamination $\mathcal{L}$. 
Let us focus on the right side of $g$ with respect to the orientation of $g$. 
Let $h_{i_1}, \cdots, h_{i_k}$ be the biinfinite leaves spiraling along $g$ from right, and $T_{i_1}, \cdots T_{i_k}$ the ideal triangles which also spiral along $g$ from right. 
Suppose that the spiralings of these leaves and triangles occur in the direction (resp. the opposite direction) of the orientation of $g$. 
Define $\overline{\sigma}_{p}^{\rho}(h_i)$ by $\sigma_p^{\rho}(h_i)$ if $h_i$ is oriented toward $g$, and by $\sigma_{n-p}^{\rho}(h_i)$ otherwise. 
We set
\begin{eqnarray*}
R_p^{\rho}(g) &=&\sum_{j=1}^{k} \overline{\sigma}_{p}^{\rho}(h_{i_j}) + \sum_{j=1}^k \sum_{q+r=n-p}\tau_{pqr}^{\rho}(T_{i_j}, v_{i_j}), \\
(\mbox{resp.  }  R_p^{\rho}(g) &=& -\sum_{j=1}^{k} \overline{\sigma}_{n-p}^{\rho}(h_{i_j}) - \sum_{j=1}^k \sum_{q+r=p}\tau_{(n-p)qr}^{\rho}(T_{i_j}, v_{i_j}), ~) 
\end{eqnarray*}
where $v_{i_j}$ is a vertex of a lift $\widetilde{T}_{i_j}$ which is a terminal (resp. starting) point of a lift of $g$. When we focus on the left side, we define $L_p^{\rho}(g)$ similarly as follows.
\begin{eqnarray*}
 L_p^{\rho}(g) &=&- \sum_{j=1}^{k} \overline{\sigma}_{p}^{\rho}(h_{i_j}) - \sum_{j=1}^k \sum_{q+r=n-p}\tau_{pqr}^{\rho}(T_{i_j}, v_{i_j}). \\
(\mbox{resp.  }  L_p^{\rho}(g) &=& \sum_{j=1}^{k} \overline{\sigma}_{n-p}^{\rho}(h_{i_j}) + \sum_{j=1}^k \sum_{q+r=p}\tau_{(n-p)qr}^{\rho}(T_{i_j}, v_{i_j}) .~) 
\end{eqnarray*}
\begin{proposition}[Bonahon-Dreyer\cite{BD2}]
For any $\rho \in \mathcal{H}_n(S)$, it holds that 
\[ l_p^{\rho}(g) = R_p^{\rho}(g) =  L_p^{\rho}(g). \]
\end{proposition}

Now we construct the parameter space of Bonahon-Dreyer.
Recall our setting: we fix
\begin{itemize} 
\item[(1)]  a hyperbolic structure on $S$ and the associated universal covering $\widetilde{S} \rightarrow S$,
\item[(2)]  a maximal geodesic oriented lamination $\mathcal{L}$ on $S$ which consists of biinfinite leaves $h_1, \cdots h_s$ and closed leaves $g_1, \dots, g_t$ with short arcs $k_1, \cdots, k_t$, and
\item[(3)] ideal triangles $T_1, \cdots, T_u$ which are defied by the ideal triangulation induced by $\mathcal{L}$.
\end{itemize}
The parameter space is given in a Euclidian space $\mathbb{R}^N$ which is described by 
\[ \mathbb{R}^N = \{ (\tau_{abc}(T_i, v_i), \cdots, \sigma_d(h_j), \cdots, \sigma_e(g_k), \cdots) \} \]
where we consider all $a,b,c,d,e,i,j$ and $k$, such that 
\begin{itemize}
\item $a,b,c,d$ and $e$ are integers with $a,b,c \geq 1$, $a+b+c =n$ and $1 \leq d,e \leq n-1$,
\item $i,j$ and $k$ are the indexes with $1\leq i \leq u$, $1\leq j \leq s$ and $1\leq k \leq t$,
\end{itemize}
and consider all vertices of a lift $\widetilde{T}_i$. 
We can count the number $N$ concretely. 
\[ N = u \cdot 3 \cdot \dfrac{(n-1)(n-2)}{2} + s \cdot (n-1) + t \cdot (n-1). \]

Let $\mathcal{P}$ be the polytope in $\mathbb{R}^N$ defined by the relations of triangle invariant given in Proposition 4.2, and the length relations $L_p=R_p>0 (p=1, \cdots, n-1)$ given in Proposition 4.5. The number of these relations is $2u(n-1)(n-2)/2 + s(n-1)$ where the first term corresponds to the relations of triangle invariants and the second corresponds the length relations. Thus, noting $u = 4g-4$ and $t=6g-6$ by the Gauss-Bonnet theorem, since the these relations are independent each other, we can compute the dimension of the convex polytope $\mathcal{P}$ as
\begin{eqnarray*}
dim \mathcal{P} 
&=& N - 2u(n-1)(n-2)/2 + s(n-1) \\
&=& u(n-1)(n-2)/2 + t(n-1) \\
&=& (4g-4)(n-1)(n-2)/2 + (6g-6)(n-1)\\
&=& (2g-2)(n^2-1).
\end{eqnarray*}

The result of Bonahon-Dreyer is as follows. 
\begin{theorem}[Bonahon-Dreyer\cite{BD2}]
The following map $\Phi_{\mathcal{L}}$ is homeomorphic onto the image.
\begin{eqnarray*}
&&\Phi_{\mathcal{L}} : {\rm Hit}_n(S) \rightarrow \mathbb{R}^N \\
&&\Phi_{\mathcal{L}}([\rho]) =  (\tau^{\rho}_{abc}(T_i, v_i), \cdots, \sigma^{\rho}_d(h_j), \cdots, \sigma^{\rho}_e(g_k), \cdots).
\end{eqnarray*}
Moreover the image of this map is the interior of the convex polytope $\mathcal{P}$
\end{theorem}

\subsection{Non-closed case}
The Bonahon-Dreyer parameterization can be extended to the case of surfaces with nonempty boundary by using Labrourie-McShane's work \cite{LM}. 
Let $S_{\partial}$ be a compact connected oriented surface with nonempty boundary, no punctures and negative Euler characteristic number. 
We denote the topological double of $S_{\partial}$ by $\widehat{S_{\partial}}$. 
Labourie-McShane studied an extension of Hitchin representations in $\mathcal{H}_n(S_{\partial})$ to $\mathcal{H}_n(\widehat{S_{\partial}})$, which is a generalization of the hyperbolic double of Fuchsian representations. 

Let $J$ be an involution of ${\rm PGL}_n(\mathbb{R})$. 
For any Hitchin representation $\rho \in \mathcal{H}_n(S_{\partial})$, Labourie-McShane constructed a unique representation $\widehat{\rho} \in \mathcal{R}_n(\widehat{S_{\partial}})$ so that $\widehat{\rho}$ is an extension of $\rho$ and the $J$-conjugate action on $\widehat{\rho}$ is compatible with the natural involution of the double $\widehat{S_{\partial}}$. 
Such a representation is called the {\it $J$-extension} of $\rho$.
In particular for the involution $J_n$ of ${\rm PGL}_n(\mathbb{R})$ defined by 
\begin{equation*}
J_n=
\begin{bmatrix}
1 & 0 & 0 &  \cdots \\
0 & -1 & 0 &  \cdots \\
0 & 0 & 1 &  \cdots \\
\vdots &\vdots & \vdots  & \ddots
\end{bmatrix},
\end{equation*}
they showed that the $J_n$-extension of a Hitchin representation is a Hitchin representation in $\mathcal{H}_n(\widehat{S_{\partial}})$. 
The $J_n$-extensions of Hitchin representations is called the {\it Hitchin doubles}. 
Note that since Hitchin doubles are representations of closed surface groups, there exist the associated flag curves for Hitchin doubles. 

\begin{theorem}[Labourie-McShane\cite{LM}]
Let $\rho \in \mathcal{H}_n(S_{\partial})$ be a Hitchin representation and $\widehat{\rho} \in \mathcal{H}_n(\widehat{S_{\partial}})$ the Hitchin double of $\rho$. 
Then the restriction to the boundary $\partial_{\infty}\widetilde{S_{\partial}}$ of the flag curve associated to $\widehat{\rho}$, $\xi_{\widehat{\rho}} : \partial_{\infty}\widetilde{\widehat{S_{\partial}}} \rightarrow {\rm Flag}(\mathbb{R}^n)$, is continuous and $\rho$-equivariant. 
\end{theorem}
We call the restriction in Theorem 4.7 the flag curve of the Hitchin representation  $\rho \in \mathcal{H}_n(S_{\partial})$. 
We remark that the proof of Theorem 4.7 implies the following lemma. 
\begin{lemma}[Labourie-McShane\cite{LM}]
The Hitchin double $\widehat{\rho_n} \in \mathcal{H}_n(S_{\partial})$ of an $n$-Fuchsian representation $\rho_n=\iota_n \circ \rho \in \mathcal{F}_n(S_{\partial})$ equals to the $n$-Fuchsian representation $\iota_n \circ \widehat{\rho} \in \mathcal{F}_n(\widehat{S_{\partial}})$ where $\widehat{\rho} \in \mathcal{F}_2(\widehat{S_{\partial}})$ is the hyperbolic double of the Fuchsian representation $\rho \in \mathcal{F}_2(S_{\partial})$.
\end{lemma}

The Bonahon-Dreyer parametrization of Hitchin components of surfaces with nonempty boundary is defined by using flag curves which are introduced above.
We can construct the parameter space similarly but we need to note invariants associated to closed geodesic leaves which are on the boundary. 
Let $\mathbb{R}^N$ denotes a Euclidian space described as follows:
\[ \mathbb{R}^N = \{ (\tau_{abc}(\widetilde{T_i}, v_i), \cdots, \sigma_d(h_j), \cdots, \sigma_e(g_k), \cdots) \} \]
where $a,b,c,d,e,i$ and $j$ satisfies the same conditions given in the closed case but $k$ is restricted so that $g_k$ is not a closed geodesic leaf on the boundary of $S_{\partial}$. 
To describe the parameter space, recall that, in the closed case, the defining polynomials are given by the relations of triangle invariants and the length relations. 
In the non-closed case, the image of the parameterization is defined by the following polynomials: (1) the relations of triangle invariants for all ideal triangles, (2) the length relations $L_p=R_p$ for all closed leaves which are not contained in the boundary, (3) the length-positivity, i.e. $R_p \mbox{~or~}L_p > 0$ for all closed leaves which are on the boundary. 
When we denote by $\mathcal{P}$ the interior of the polytope defined by the relations (1), (2) and (3), Theorem 4.6 can be shown similarly.


\section{Computation of Fuchsian loci}
\subsection{A parametrization of Fuchsian representations}
We determine the Fuchsian locus of a pair of pants by using the Bonahon-Dreyer coordinate. 
Let $P$ be a pair of pants with the maximal geodesic lamination $\mathcal{L}=\{h_{AB}, h_{BC}, h_{CA}, A, B, C \}$ which is described in Figure 5. 
This lamination induces an ideal triangulation of $P$ and we denote by $T_0$ and $T_1$ two triangles given by the triangulation.  
Orientations of leaves are fixed as Figure 5. 
Let $a,b$ and $c$ be the homotopy classes of the boundary  components $A,B$ and $C$ of $P$.
Then the fundamental group $\pi_1(P)$ is generated by $a,b$ and $c$, and has the presentation
\[ \pi_1(P) = <a,b,c~|~abc=1>. \]

\begin{figure}[htbp]
\begin{center}
\includegraphics[width = 3cm,height=3cm]{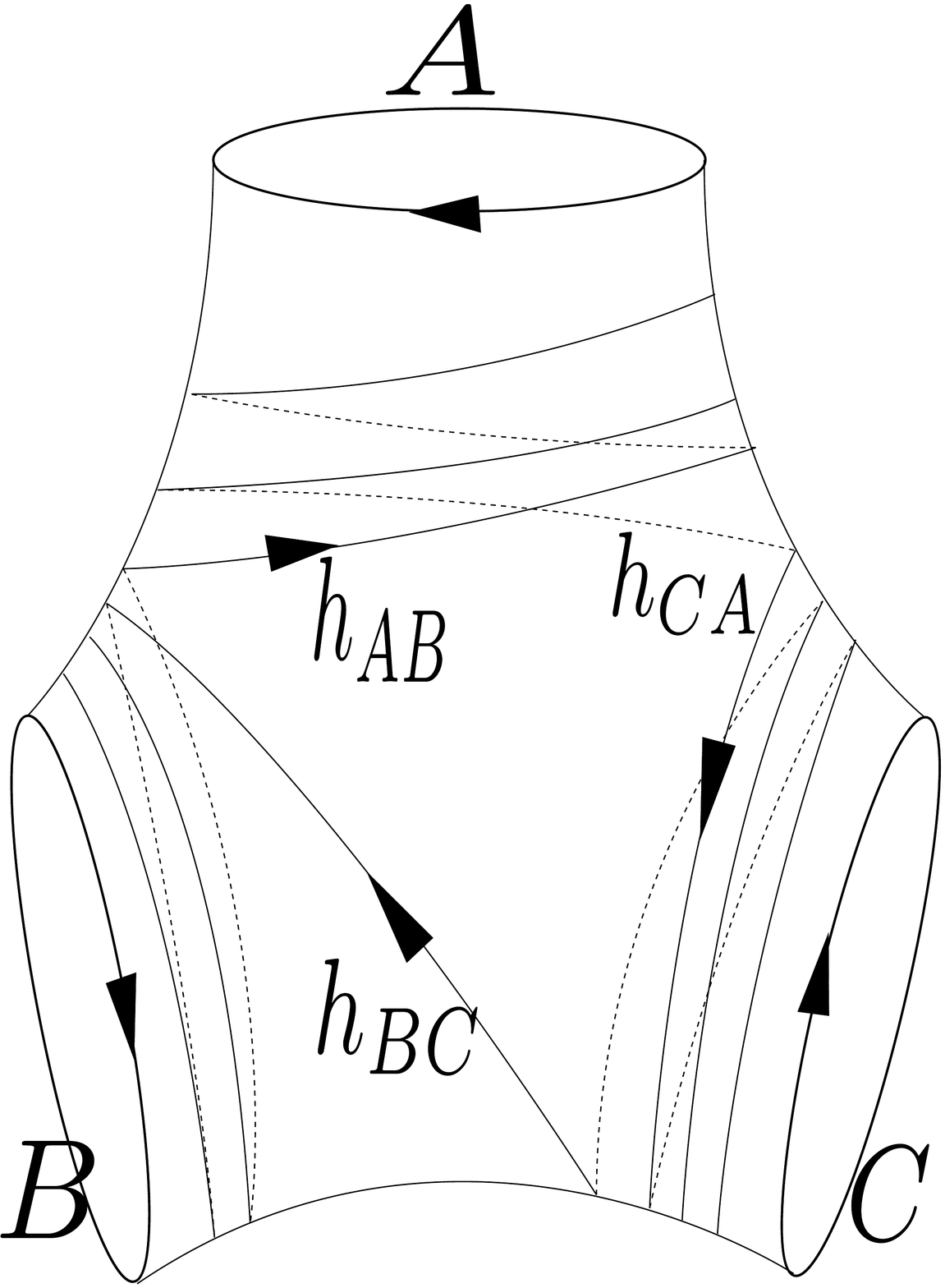}
\caption{The lamination $\mathcal{L}$ on $P$.}
\label{figure 7}
\end{center}
\end{figure}

At first we parametrize the Fuchsian representations $\rho : \pi_1(P) \rightarrow {\rm PSL}_2(\mathbb{R})$ by hyperbolic structures of $P$.
As we saw in Proposition 2.1 the Teichm\"uller space $\mathscr{T}(P)$ is parametrized by the hyperbolic lengths of the boundary components of $P$, and hence the Fuchsian representations are also uniquely determined by the lengths up to conjugacy. 

\begin{proposition}
Let ${\bf m}=(l_A, l_B, l_C)$ be a triple of the hyperbolic lengths of the boundary components $A,B$ and $C$. 
Then we can take a representative $\rho$ in the conjugacy class of the Fuchsian representation associated to the data ${\bf m}$ such that the developing map ${\bf dev}_{\rho}$ associated to $\rho$ is described as Figure 6 and 7, and the attracting point of the axis of $\rho(b)$ equals to $0$ in $\partial_{\infty}\mathbb{H}^2$. 
In particular, the biinfinite leaves $h_{AB}, h_{BC}$ and $h_{CA}$ can lift to the geodesics $\widetilde{h}_{AB}, \widetilde{h}_{BC}$ and $\widetilde{h}_{CA}$ in Figure 6 and 7, whose terminal points are $\infty, 1$ and $0$ in $\partial_{\infty}\mathbb{H}^2$ respectively. 
Moreover we can write such a representative $\rho$ concretely as follows.
\[
\rho(a) =
\begin{bmatrix}
\alpha & \alpha \beta \gamma + \alpha^{-1} \\
0 & \alpha^{-1}
\end{bmatrix},~~
\rho(b) =
\begin{bmatrix}
\gamma & 0 \\
-\beta^{-1}-\gamma^{-1} & \gamma^{-1}
\end{bmatrix},
\]
where $\alpha, \beta, \gamma : \mathbb{R}_{>0}^3 \rightarrow \mathbb{R}_{>0}$ are defined by
\[
\alpha(l_A, l_B, l_C) = e^{l_A/2},~~
\beta(l_A, l_B, l_C)  = e^{(l_C - l_A)/2},~~
\gamma(l_A, l_B, l_C) =  e^{-l_B/2}
\]
with the conditions $\alpha >1, 1> \gamma >0, \beta >0$.  
\end{proposition}

\begin{figure}[htbp]
\begin{minipage}{0.45\hsize}
\begin{center}
\includegraphics[width = 45mm]{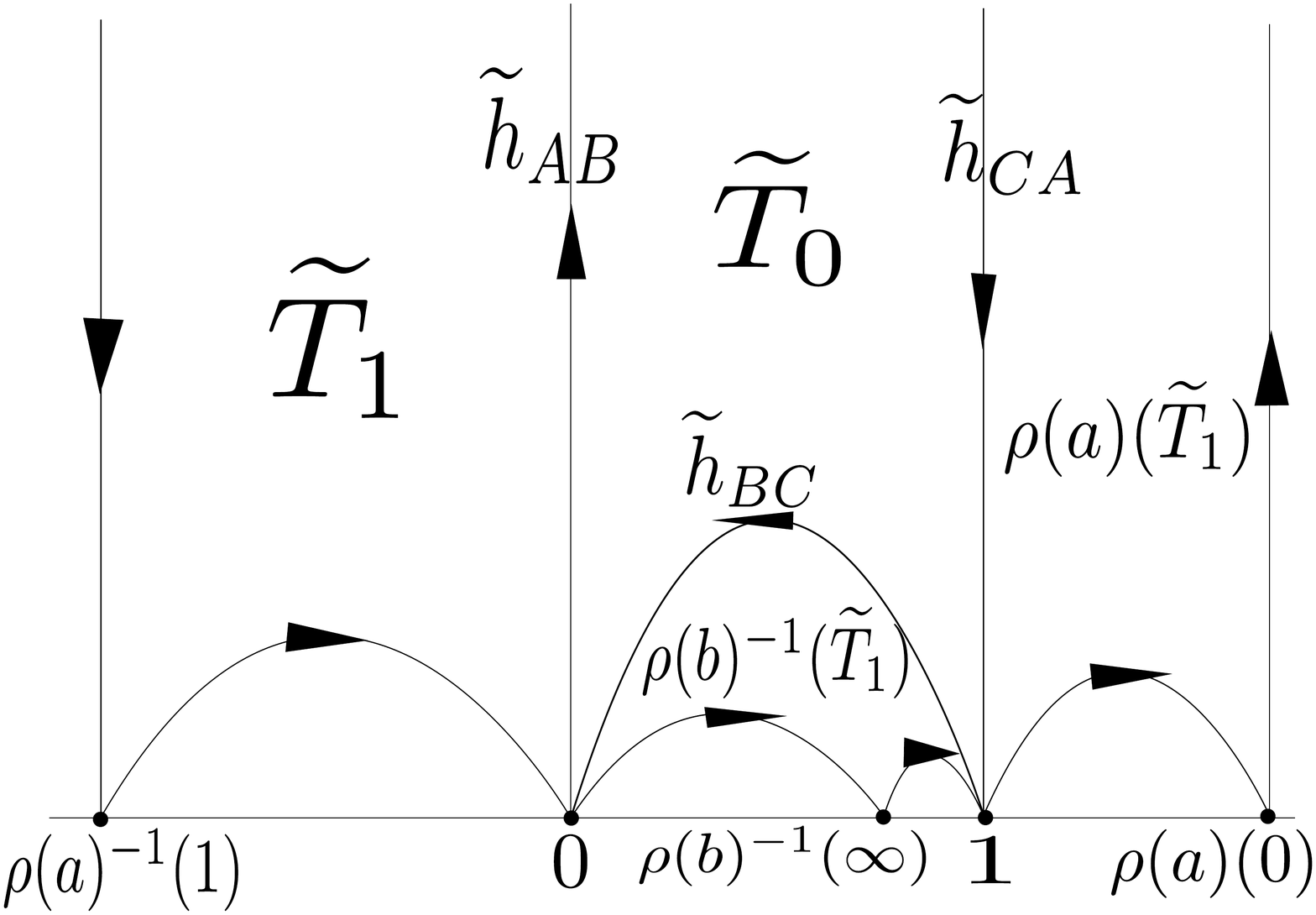}
\end{center}
\caption{Upper half plane model.}
\end{minipage}
\begin{minipage}{0.45\hsize}
\begin{center}
\includegraphics[width=45mm, height=45mm]{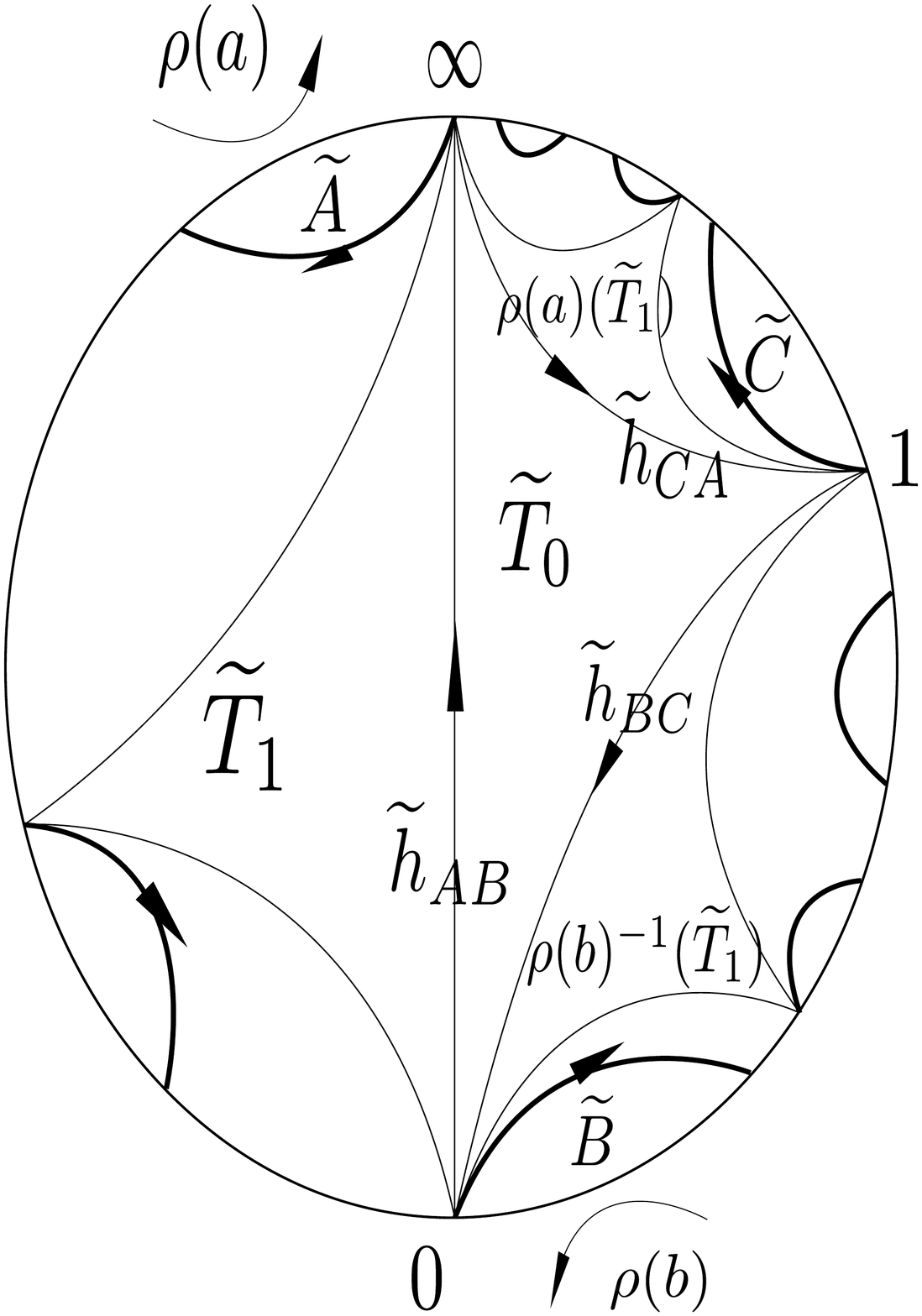}
\end{center}
\caption{Poincar\'e disk model.}
\end{minipage}
\end{figure}
\begin{flushleft}
\underline{{\it Proof.}}
\end{flushleft}
The first assertion follows by the normalization of Fuchsian groups. See \cite{IT}. 
We can set the Fuchsian representation $\rho$ which satisfies the condition of the fixed points of $\rho(a)$ and $\rho(b)$ by
\[
\rho(a) =
\begin{bmatrix}
\alpha & \beta \\
0 & \alpha^{-1}
\end{bmatrix},~~
\rho(b) =
\begin{bmatrix}
\gamma & 0 \\
\delta & \gamma^{-1}
\end{bmatrix}.
\]
First we note that $\alpha$ and $\beta$ satisfies $-\alpha^{-1} + \beta \neq 0$, since if $\alpha^{-1} = \beta$ then 
\[
\rho(c) = \rho(b)^{-1} \rho(a)^{-1} =
\begin{bmatrix}
\alpha^{-1}\gamma^{-1} & -\alpha^{-1} \gamma^{-1} \\
-\alpha^{-1} \delta & \alpha^{-1} \delta + \alpha \gamma
\end{bmatrix},
\]
and 
\[ \rho(c)(1) = \dfrac{0}{0 + \alpha\gamma} = 0 \]
where $\alpha\gamma \neq 0$.
This contradicts the assumption that $\rho(c)$ fixes 1. 
Therefore we may assume $-\alpha^{-1} + \beta \neq 0$
We compute the parameter $\delta$ so that $\rho(c) = \rho(b)^{-1} \rho(a)^{-1}$ fixes $1 \in \partial_{\infty}\mathbb{H}^2$. Since
\[
\rho(c) = \rho(b)^{-1} \rho(a)^{-1} =
\begin{bmatrix}
\alpha^{-1}\gamma^{-1} & -\beta \gamma^{-1} \\
-\alpha^{-1} \delta & \beta \delta + \alpha \gamma
\end{bmatrix},
\]
the condition of the fixed points of $\rho(c)$ implies that
\[ \rho(c)(1) = \dfrac{\alpha^{-1} \gamma^{-1} - \beta \gamma^{-1}}{-\alpha^{-1}\delta + \beta \delta + \alpha \gamma} = 1. \] 
Thus we obtain
\[ \delta = -\gamma^{-1} + \dfrac{ - \alpha \gamma}{-\alpha^{-1} + \beta} .\]
Replacing the parameter $\beta$ with a parameter $\beta'$ satisfying $\beta = \alpha \gamma \beta'  + \alpha^{-1} $, we deform the equation above as 
\[ \delta = -\gamma^{-1} - \beta'^{-1} .\]
We denote $\beta'$ by $\beta$ again. The representation $\rho$ is given by the following form.
\[
\rho(a) =
\begin{bmatrix}
\alpha & \alpha \beta \gamma + \alpha^{-1} \\
0 & \alpha^{-1}
\end{bmatrix},~~
\rho(b) =
\begin{bmatrix}
\gamma & 0 \\
-\beta^{-1}-\gamma^{-1} & \gamma^{-1}
\end{bmatrix}.
\]
Note that $\rho(a)$ and $\rho(b)$ are elements in ${\rm PSL}_2(\mathbb{R})$, so we can assume $\alpha, \gamma >0$ by the multiple of $-1$. Since the attracting point of $\rho(b)$ is $0 \in \partial_{\infty}\mathbb{H}^2$, it holds that $\gamma < 1$. 

Next we consider inequalities among parameters $\alpha, \beta$ and $\gamma$ which are given by positional relations of fixed points at infinity. 
We can check the following easily:
\begin{eqnarray*}
{\rm Fix}(\rho(a)) &=& \{ \infty, \dfrac{\alpha^2\beta\gamma+1}{1-\alpha^2} \} ,\\
{\rm Fix}(\rho(b)) &=& \{ 0, \dfrac{\gamma - \gamma^{-1}}{-\beta^{-1} - \gamma^{-1}} \} ,\\
{\rm Fix}(\rho(c)) &=& \{ 1, \dfrac{\alpha \beta + \alpha^{-1}\gamma^{-1}}{\alpha^{-1}\gamma^{-1} + \alpha{-1}\beta^{-1}} \}.
\end{eqnarray*}
By Figure 6, the following inequalities hold:
\begin{eqnarray}
&& \dfrac{\alpha^2\beta\gamma+1}{1-\alpha^2} < 0 ,\\
&& 0 < \dfrac{\gamma - \gamma^{-1}}{-\beta^{-1} - \gamma^{-1}} < 1 ,\\
&& 1< \dfrac{\alpha \beta + \alpha^{-1}\gamma^{-1}}{\alpha^{-1}\gamma^{-1} + \alpha{-1}\beta^{-1} }.
\end{eqnarray}
Noting that $1> \gamma >0$ which implies that $\gamma < \gamma^{-1}$, these inequalities are deformed as follows:
\begin{eqnarray}
(2)  &\Leftrightarrow& \beta^{-1} + \gamma^{-1}  >0, \beta^{-1} + \gamma >0 ,\\
(3)+(4) &\Leftrightarrow& \alpha^2 \beta > \beta^{-1}, \\
(1)+(4) &\Leftrightarrow& \alpha > 1 .
\end{eqnarray}
We can deduce $\beta > 0$ from these inequalities since if $\beta$ is negative, 
\begin{eqnarray*}
(4)  &\Leftrightarrow& \beta + \gamma < 0, \beta < - \dfrac{1}{\gamma}, \\
(5) &\Leftrightarrow& -\dfrac{1}{\alpha}< \beta < \dfrac{1}{\alpha}
\end{eqnarray*}
and, by $-1/\gamma < -1 / \alpha$, a contradiction occurs. 

In general, the hyperbolic length $l$ of a simple closed geodesic which is covered by the axis of a hyperbolic isometry $M \in {\rm PSL}_2(\mathbb{R})$ is given by the following formula:
\[ |tr(M) | = 2\cosh (l/2) .\]
Under the conditions above of $\alpha, \beta$ and $\gamma$, the data of the hyperbolic lengths of $A$ and $B$ detect $\alpha$ and $\gamma$.
\[ \alpha = \exp (l_A/2),~~\gamma = \exp (-l_B/2) .\]
We consider an equation 
\[ |tr(\rho(c))| = 2\cosh (l_C/2) \]
which implies that 
\begin{eqnarray*}
&& |-\alpha \beta - \alpha^{-1} \beta^{-1} | = 2 (\cosh(l_C/2)) \\
&\Leftrightarrow&
\alpha \beta + \alpha^{-1} \beta^{-1} = 2\cosh(l_C/2)~~~~(\alpha, \beta >0) \\
&\Leftrightarrow&
\beta = \dfrac{\cosh(l_C/2) \pm \sinh(l_C/2)}{\alpha}.
\end{eqnarray*} 
The inequality (5) gives us the following condition
\[ (5) \Leftrightarrow \alpha^2 \beta^2 >1 \]
and then $\beta$ is uniquely determined by
\[ \beta = \dfrac{\cosh(l_C/2) + \sinh(l_C/2)}{\alpha} .\]
\begin{flushright}
$\Box$
\end{flushright}
\subsection{A parametrization of Fuchsian loci}
Let $\rho \in \mathcal{F}_2(P)$ be a Fuchsian representation of the form in Proposition 5.1, and $\rho_n = \iota_n \circ \rho$ be the induced $n$-Fuchsian representation in $\mathcal{F}_n(P)$. 
Since $P$ has the nonempty boundary $A,B$ and $C$ and the double of $P$ is the closed surface $\Sigma_2$ of genus 2, we must consider the Hitchin double $\widehat{\rho_n} \in \mathcal{H}_n(\Sigma_2)$ of $\rho_n$ to construct the flag curve $\xi_{\rho_n}$ of $\rho_n$ which is the restriction of the flag curve $\xi_{\widehat{\rho_n}}$ of $\widehat{\rho_n}$. 
The Hitchin double of $n$-Fuchsian representation is an $n$-Fuchsian by Lemma 4.8, in particular it holds that $\widehat{\rho_n}=\iota_n \circ \widehat{\rho}$ where $\widehat{\rho} \in \mathcal{F}_2(\Sigma_2)$ is the $J_2$-extension of $\rho$, that is, the hyperbolic double of $\rho$. 
Therefore it is enough to consider the flag curve of the $n$-Fuchsian representation. 

In the case of $n$-Fuchsian representations, we can explicitly describe the associated flag curves.
Following \cite{Gu}, we describe the flag curves of $n$-Fuchsian representations $\widehat{\rho_n}=\iota_n \circ \widehat{\rho}$. The developing map associated to the holonomy representation $\widehat{\rho}$ induces an identification of the universal covering $\widetilde{\Sigma_2}$ with $\mathbb{H}^2$. 
Note that this identification is extended to the identification between the visual boundaries of $\widetilde{\Sigma_2}$ and $\mathbb{H}^2$.
We identify $\mathbb{R}^n$ with the space of homogeneous polynomials $F(X,Y)$ of degree $n-1$ as vector spaces.
Set
\[ V=Span_{\mathbb{R}}<X^{n-1}, X^{n-2}Y, \cdots, Y^{n-1}> .\]
For $i=1, \cdots, n$ and $r \in \partial_{\infty}\mathbb{H}^2$, $W^{(i)}(r)$ denotes the subspace of $V$ consisting of polynomials which can be divided by $(rX+Y)^{n-i}$ if $r \neq \infty$ and by $X^{n-i}$ if $r=\infty$. 
Set $W^{(0)}(r) = 0$ for any $r$. 
Consider the map 
\[
\xi: \partial_{\infty}\mathbb{H}^2 \rightarrow {\rm Flag}(\mathbb{R}^n)~:~  r \mapsto \{W^{(i)}(r) \}_{i=0, \cdots, n}.
\]
It is known that the composition $\xi_{\widehat{\rho_n}} : \partial_{\infty}\widetilde{\Sigma_2} = \partial_{\infty}\mathbb{H}^2 \rightarrow {\rm Flag}(\mathbb{R}^n)$ is the flag curve of the Hitchin representation $\widehat{\rho_n}$.

By Theorem 4.7, the flag curve associated to $\rho_n$ is the restriction of $\xi_{\widehat{\rho_n}}$ to $\partial_{\infty}\widetilde{P}$. 
Let ${\bf dev}_{\widehat{\rho}}$ be the developing map defined by $\widehat{\rho}$, and $i:\widetilde{P} \rightarrow \widetilde{\Sigma_2}$ the natural inclusion. 
Then it holds that ${\bf dev}_{\rho} = {\bf dev}_{\widehat{\rho}} \circ i$ since $\widehat{\rho}$ is the hyperbolic double of $\rho$. 
Thus $\xi_{\rho_n}$ is given by 
\[ \xi_{\rho_n}: \partial_{\infty}P \xrightarrow{{\bf dev}_{\rho}} \partial_{\infty}\mathbb{H}^2 \xrightarrow{\xi} {\rm Flag}(\mathbb{R}^n) \]


Recalling the definition of the Bonahon-Dreyer parameterization, it is enough to compute the shearing invariants $\sigma_p^{\rho_n}(h_{AB}), \sigma_p^{\rho_n}(h_{BC}), \sigma_p^{\rho_n}(h_{CA})~(1 \leq p \leq n-1)$ and the triangle invariants $\tau_{pqr}^{\rho_n}(T_0, \infty), \tau_{pqr}^{\rho_n}(T_1, \infty)~(p,q,r \geq 1~s. t.~p+q+r=n)$. 
These invariants detect the image of the Fuchsian locus of a pair of pants by the parametrization $\Phi_{\mathcal{L}}$. 
In the following computation, we use a useful formula below.

\begin{lemma}
Let $V$ be an n-dimensional vector space with a basis $\{ b_1, \cdots, b_n \}$, and $v_1, \cdots, v_n$ be arbitrary vectors in $V$. 
When $v_j = \sum_{i=1}^n v_{ij}b_j$, it holds that
\[ v_1 \wedge \cdots \wedge v_n=Det((v_{ij})_{1 \leq i,j \leq n}) \cdot b_1 \wedge \cdots \wedge b_n .\]  
\end{lemma}
\begin{flushright}
$\Box$
\end{flushright}
Now we may fix an identification $\bigwedge^{(n)} V \rightarrow \mathbb{R}$ so that $X^{n-1} \wedge X^{n-2}Y \wedge \cdots \wedge Y^{n-1}$ is identified with $1 \in \mathbb{R}$.
Moreover we extend the definition of binomial coefficients formally as follows.
\[
\binom{n}{p} =
\begin{cases}
\dfrac{n!}{p!(n-p)!}  & (0 \leq p \leq n) \\
0 & (otherwise).
\end{cases}
\]
Fix the basis $b_1=X^{n-1}, b_2 = X^{n-2}Y, \cdots, b_n=Y^{n-1}$ of the $n$-dimensional vector space $V=\mathbb{R}[X^{n-1}, X^{n-2}Y, \cdots, Y^{n-1}]$. 
\begin{flushleft}
\fbox{1.Triangle invariant $\tau_{pqr}^{\rho_n}(T_0, \infty)$}
\end{flushleft}
Consider a lifting $\widetilde{T_0}$ of the ideal triangle $T_0$ whose vertices are $\infty$, $0$ and $1$, and denote the corresponding flags  $E=\xi_{\rho_n}(\infty), F=\xi_{\rho_n}(1)$ and $G=\xi_{\rho_n}(0)$. 
To detect the triangle invariant of the ideal triangle $\widetilde{T}_0$ we observe the value $X(p, q, r) = e^{(p)} \wedge f^{(q)} \wedge g^{(r)}$ (Definition 3.2). By the definition of $\xi_{\rho_n}$, the flags $E^{(p)}$, $F^{(q)}$ and $G^{(r)}$ are described as follows:
{\small
\begin{eqnarray*}
E^{(p)} &=& {\rm Span}_{\mathbb{R}} <X^{n-1}, X^{n-2}Y, \cdots, X^{n-p}Y^{p-1}> ,\\
F^{(q)} &=& {\rm Span}_{\mathbb{R}} <(X+Y)^{n-q}X^{q-1}, (X+Y)^{n-q}X^{q-2}Y, \cdots, (X+Y)^{n-q}Y^{q-1}> ,\\
G^{(r)} &=& {\rm Span}_{\mathbb{R}} <Y^{n-1}, XY^{n-2}, \cdots, X^{r-1}Y^{n-r}> .
\end{eqnarray*}
}
Therefore we may choose nonzero elements as below:
{\small
\begin{eqnarray*}
e^{(p)} &=& X^{n-1}\wedge X^{n-2}Y \wedge \cdots \wedge X^{n-p}Y^{p-1} ,\\
f^{(q)} &=& (X+Y)^{n-q}X^{q-1} \wedge (X+Y)^{n-q}X^{q-2}Y \wedge \cdots \wedge (X+Y)^{n-q}Y^{q-1}  ,\\
g^{(r)} &=& Y^{n-1} \wedge XY^{n-2} \wedge \cdots \wedge X^{r-1}Y^{n-r} .
\end{eqnarray*}
}

To apply Lemma 5.2, we write the vectors $e^{(p)},f^{(q)},g^{(r)} $ by their entries with respect to the basis $\{b_i\}$. 
Note that
{\small
\begin{eqnarray*}
&&(X+Y)^{n-q}X^{q-k}Y^{k-1} \\
&=& 
(\binom{n-q}{0} X^{n-q} +  \binom{n-q}{1} X^{n-q-1}Y + \cdots +  \binom{n-q}{n-q} Y^{n-q}) X^{q-k}Y^{k-1}\\
&=&
\binom{n-q}{0}b_k + \binom{n-q}{1}b_{k+1} + \cdots +  \binom{n-q}{n-q}b_{n-q+k}.
 \end{eqnarray*}
}
Setting
\[
A_{pqr}(T_0)=
\begin{pmatrix}
\binom{p+r}{0}      && \binom{p+r}{-1}      && \cdots && \binom{p+r}{-q+1} \\
\binom{p+r}{1}      && \binom{p+r}{0}        && \cdots && \binom{p+r}{-q+2}	\\
\vdots &&  \vdots && \vdots                 && \vdots   			\\
\binom{p+r}{p+r}   && \binom{p+r}{p+r-1} && \cdots && \binom{p+r}{0}		\\
\binom{p+r}{p+r+1}&& \binom{p+r}{p+r}    && \cdots && \binom{p+r}{1}		\\
\vdots                 && \vdots                 && \vdots  && \vdots 			\\
\binom{p+r}{n-1}  && \binom{p+r}{n-2}    && \cdots && \binom{p+r}{p+r}	\\
\end{pmatrix},
\]
and using Lemma 5.2, we compute $X(p,q,r)$: for $q \neq 0$
\begin{eqnarray*}
X(p,q,r)
&=& 
e^{(p)} \wedge f^{(q)} \wedge g^{(r)}\\
&=&
\begin{vmatrix}
{\rm Id}_p &  & \mbox{\LARGE 0} \\
  & \mbox{\huge $A_{pqr}(T_0)$ } & \\ 
\mbox{\LARGE 0} &  & {\rm Id}_r
\end{vmatrix}\\
&=&
\begin{vmatrix}
\binom{p+r}{p} & \cdots & \binom{p+r}{p-q+1} \\
\vdots & \vdots & \vdots \\
\binom{p+r}{p+q-1} & \cdots & \binom{p+r}{p} \\
\end{vmatrix}
\end{eqnarray*}
and $X(p,0,r) = 1$ for all $p,r$.

\begin{flushleft}
\fbox{2.Triangle invariant $\tau_{pqr}^{\rho_n}(T_1, \infty)$}
\end{flushleft}
This computation is similar to the case of the triangle invariant of $T_0$. 
We take a lifting $\widetilde{T}_1$ of $T_1$ so that the vartices are $\infty, 0$ and $(\rho(a))^{-1}(1)=-\beta \gamma$. 
We set as $E= \xi_{\rho_n}(\infty), F=\xi_{\rho_n}(0)$ as $G = \xi_{\rho_n}(-\beta \gamma)$. 
Then these flags are described by 
{\small
\begin{eqnarray*}
E^{(p)} &=& {\rm Span}_{\mathbb{R}} <X^{n-1}, X^{n-2}Y, \cdots, X^{n-p}Y^{p-1}>, \\
F^{(q)} &=& {\rm Span}_{\mathbb{R}} <Y^{n-1}, XY^{n-2}, \cdots, X^{q-1}Y^{n-q}> , \\
G^{(r)} &=& {\rm Span}_{\mathbb{R}} <(-\beta\gamma X+Y)^{n-q}X^{q-1}, (-\beta \gamma X+Y)^{n-q}X^{q-2}Y, \cdots \\
&& \hspace{10pt} \cdots, (-\beta \gamma X+Y)^{n-q}Y^{q-1}> ,
\end{eqnarray*}
}
and nonzero elements in $X(p,q,r)$ are chosen as  
{\small
\begin{eqnarray*}
e^{(p)} &=& X^{n-1}\wedge X^{n-2}Y \wedge \cdots \wedge X^{n-p}Y^{p-1} ,\\
f^{(q)} &=& Y^{n-1} \wedge XY^{n-2} \wedge \cdots \wedge X^{q-1}Y^{n-q}  ,\\
g^{(r)} &=& (-\beta\gamma X+Y)^{n-r}X^{r-1} \wedge (-\beta \gamma X+Y)^{n-r}X^{r-2}Y \wedge \cdots \wedge (-\beta \gamma X+Y)^{n-r}Y^{r-1} .
\end{eqnarray*}
}
Set
\begin{equation*}
A_{pqr}(T_1)=
\begin{pmatrix}
\binom{p+q}{0}(-\beta  \gamma)^{p+q}     && \binom{p+q}{-1}(-\beta  \gamma)^{p+q+1}      && \cdots && \binom{p+q}{-r+1}(-\beta  \gamma)^{n-1} \\
\binom{p+q}{1}(-\beta  \gamma)^{p+q-1}      && \binom{p+q}{0}(-\beta  \gamma)^{p+q}        && \cdots && \binom{p+q}{-r+2}(-\beta  \gamma)^{n-2}	\\
\vdots &&  \vdots && \vdots                 && \vdots   			\\
\binom{p+q}{p+q}(-\beta  \gamma)^{0}   && \binom{p+q}{p+q-1}(-\beta  \gamma)^{1} && \cdots && \binom{p+q}{0}	(-\beta  \gamma)^{p+q}	\\
\binom{p+q}{p+q+1} (-\beta  \gamma)^{-1} && \binom{p+q}{p+q}(-\beta  \gamma)^{0}    && \cdots && \binom{p+q}{1}	(-\beta  \gamma)^{p+q-1}	\\
\vdots                 && \vdots                 && \vdots  && \vdots 			\\
\binom{p+q}{n-1} (-\beta  \gamma)^{-r+1} && \binom{p+q}{n-2}(-\beta  \gamma)^{-r+2}    && \cdots && \binom{p+q}{p+q}(-\beta  \gamma)^{0}	\\
\end{pmatrix}.
\end{equation*}
Therefore, for $r \neq 0$
\begin{eqnarray*}
X(p,q,r)
&=& 
e^{(p)} \wedge f^{(q)} \wedge g^{(r)}\\
&=& 
\begin{vmatrix}
{\rm Id}_p &  \mbox{\LARGE 0} & \\
 & & \mbox{\huge $A_{pqr}(T_1)$ }  \\ 
\mbox{\LARGE 0} &  {\rm Id}_q & 
\end{vmatrix}\\
&=& (-1)^{q(r+1)}
\begin{vmatrix}
\binom{p+q}{p}(-\beta  \gamma)^{q} & \cdots & \binom{p+q}{p-r+1}(-\beta  \gamma)^{q+r-1} \\
\vdots & \vdots & \vdots \\
\binom{p+q}{p+r-1}(-\beta  \gamma)^{q-r+1} & \cdots & \binom{p+q}{p}(-\beta  \gamma)^{q} \\
\end{vmatrix},
\end{eqnarray*}
and $ X(p,q,0) = (-1)^q$ for all $p,q$. 
\begin{flushleft}
\fbox{3.Shearing invariant $\sigma_p^{\rho_n}(h_{AB})$}
\end{flushleft}

We take a lifting of $h_{AB}$ so that the starting and terminal point of $h_{AB}$ are $0$ and $\infty$ respectively, and consider the rectangle whose vertices are $\infty, -\beta \gamma, 0$ and $1$ where the point $-\beta\gamma$ is the image of the covering translation $ (\rho(a))^{-1}(1) = -\beta \gamma$. 
Then the shearing invariant of $h_{AB}$ is defined by
\[ \sigma_p^{\rho_n}(h_{AB}) = \log D_p (\xi_{\rho_n}(\infty), \xi_{\rho_n}(0), \xi_{\rho_n}(-\beta \gamma), \xi_{\rho_n}(1)) . \]
We set $E=\xi_{\rho_n}(\infty), F=\xi_{\rho_n}(0), G=\xi_{\rho_n}(-\beta \gamma)$ and $G'=\xi_{\rho_n}(1)$.
These flags are concretely given by
{\small
\begin{eqnarray*}
E^{(p)} &=& {\rm Span}_{\mathbb{R}} <X^{n-1}, X^{n-2}Y, \cdots, X^{n-p}Y^{p-1}> ,\\
F^{(p)} &=& {\rm Span}_{\mathbb{R}} <Y^{n-1}, XY^{n-2}, \cdots, X^{p-1}Y^{n-p}> , \\
G^{(p)} &=& {\rm Span}_{\mathbb{R}} <(-\beta\gamma X+Y)^{n-p}X^{p-1}, (-\beta \gamma X+Y)^{n-p}X^{p-2}Y, \cdots \\
&& \hspace{10pt} \cdots, (-\beta \gamma X+Y)^{n-p}Y^{p-1}> ,\\
G'^{(p)} &=& {\rm Span}_{\mathbb{R}} <(X+Y)^{n-p}X^{p-1}, (X+Y)^{n-p}X^{p-2}Y, \cdots \\
&& \hspace{10pt} \cdots, (X+Y)^{n-p}Y^{p-1}>.
\end{eqnarray*}
}
We may chose nonzero elements as
{\small
\begin{eqnarray*}
e^{(p)} &=& X^{n-1} \wedge X^{n-2}Y \wedge \cdots \wedge X^{n-p}Y^{p-1} ,\\
f^{(p)} &=& Y^{n-1} \wedge XY^{n-2} \wedge \cdots \wedge X^{p-1}Y^{n-p}  ,\\
g^{(p)} &=& (-\beta\gamma X+Y)^{n-p}X^{p-1} \wedge (-\beta \gamma X+Y)^{n-p}X^{p-2}Y \wedge \cdots \\
&& \hspace{10pt} \cdots \wedge (-\beta \gamma X+Y)^{n-p}Y^{p-1} ,\\
g'^{(p)} &=& (X+Y)^{n-p}X^{p-1} \wedge (X+Y)^{n-p}X^{p-2}Y \wedge \cdots \\
&& \hspace{10pt} \cdots \wedge (X+Y)^{n-p}Y^{p-1}.
\end{eqnarray*}
}
Under the basis $b_1, \cdots, b_n$, the vectors $g^{(1)}$ and $g'^{(1)}$ are presented by their entries as follows:
{\small
\begin{eqnarray*}
g^{(1)} &=& \binom{n-1}{0}(-\beta \gamma)^{n-1} b_1 + \binom{n-1}{1}(-\beta \gamma)^{n-2} b_2 + \cdots + \binom{n-1}{n-1}(-\beta \gamma)^{0} b_n ,\\
g'^{(1)} &=&  \binom{n-1}{0} b_1 + \binom{n-1}{1} b_2 + \cdots + \binom{n-1}{n-1} b_n .
\end{eqnarray*}
}
Then we compute $Y(p)$ in Definition 3.4 as follows.
{\small
\begin{eqnarray*}
Y(p) &=& e^{(p)} \wedge f^{(n-p-1)} \wedge g^{(1)} \\
&=&
\begin{vmatrix}
{\rm Id}_p &  \mbox{\LARGE 0} & \binom{n-1}{0}(-\beta \gamma)^{n-1} \\
 & &   \binom{n-1}{1}(-\beta \gamma)^{n-2} \\
 & & \vdots \\ 
\mbox{\LARGE 0} &  {\rm Id}_{n-p-1} & \binom{n-1}{n-1}(-\beta \gamma)^{0}
\end{vmatrix}\\
&=&
(-1)^{n-p-1}\binom{n-1}{p}(-\beta \gamma)^{n-p-1} \\
&=&
\binom{n-1}{p}(\beta \gamma)^{n-p-1}.
\end{eqnarray*}
}

Similarly $Y'(p)$ in Definition 3.4 is computed as follows;
{\small
\begin{eqnarray*}
Y'(p) &=& e^{(p)} \wedge f^{(n-p-1)} \wedge g'^{(1)} \\
&=&
\begin{vmatrix}
{\rm Id}_p &  \mbox{\LARGE 0} & \binom{n-1}{0} \\
 & &   \binom{n-1}{1} \\
 & & \vdots \\ 
\mbox{\LARGE 0} &  {\rm Id}_{n-p-1} & \binom{n-1}{n-1}
\end{vmatrix}\\
&=&
(-1)^{n-p-1}\binom{n-1}{p} .\\
\end{eqnarray*}
}
\begin{flushleft}
\fbox{4.Shearing invariant $\sigma_p^{\rho_n}(h_{BC})$}
\end{flushleft}

This invariant is similarly computed as the invariant $\sigma_p^{\rho_n}(h_{AB})$.
We lift $h_{BC}$ to $\widetilde{h}_{BC}$ so that the starting and terminal points are $1$ and $0$, and consider the rectangle whose vertices are $1, \infty, 0$ and $\beta/(\beta+\gamma)$ where $(\rho(b))^{-1}(\infty)  = \beta/(\beta+\gamma)$. 
We set $E=\xi_{\rho_n}(0), F=\xi_{\rho_n}(1), G=\xi_{\rho_n}(\beta/(\beta+\gamma))$ and $G'=\xi_{\rho_n}(\infty)$. 
These flags are described by
{\small
\begin{eqnarray*}
E^{(p)} &=& {\rm Span}_{\mathbb{R}} <Y^{n-1}, XY^{n-2}, \cdots, X^{p-1}Y^{n-p}> ,\\
F^{(p)} &=&  {\rm Span}_{\mathbb{R}} <(X+Y)^{n-p}X^{p-1}, (X+Y)^{n-p}X^{p-2}Y, \cdots, (X+Y)^{n-p}Y^{p-1}> ,\\
G^{(p)} &=& {\rm Span}_{\mathbb{R}} <(\dfrac{\beta}{\beta+\gamma} X+Y)^{n-p}X^{p-1}, (\dfrac{\beta}{\beta+\gamma} X+Y)^{n-p}X^{p-2}Y, \cdots \\
&& \hspace{10pt} \cdots, (\dfrac{\beta}{\beta+\gamma} X+Y)^{n-p}Y^{p-1}>, \\
G'^{(p)} &=&   {\rm Span}_{\mathbb{R}} <X^{n-1}, X^{n-2}Y, \cdots, X^{n-p}Y^{p-1}>,
\end{eqnarray*}
}
We choose nonzero elements as
{\small
\begin{eqnarray*}
e^{(p)} &=& Y^{n-1} \wedge XY^{n-2} \wedge \cdots \wedge X^{p-1}Y^{n-p} ,\\
f^{(p)} &=& (X+Y)^{n-p}X^{p-1} \wedge (X+Y)^{n-p}X^{p-2}Y \wedge \cdots \wedge (X+Y)^{n-p}Y^{p-1} ,\\
g^{(p)} &=& Y^{n-1}  (\dfrac{\beta}{\beta+\gamma} X+Y)^{n-p}X^{p-1} \wedge (\dfrac{\beta}{\beta+\gamma} X+Y)^{n-p}X^{p-2}Y \wedge \cdots \\
&& \hspace{10pt} \cdots \wedge (\dfrac{\beta}{\beta+\gamma} X+Y)^{n-p}Y^{p-1} ,\\
g'^{(p)} &=&  X^{n-1} \wedge X^{n-2}Y \wedge \cdots \wedge X^{n-p}Y^{p-1} .\\
\end{eqnarray*}
}
The vectors $g^{(1)}$ and $g'^{(1)}$ are presented by their entries as follows:
{\small
\begin{eqnarray*}
g^{(1)} &=& \binom{n-1}{0}(\dfrac{\beta}{\beta+\gamma})^{n-1} b_1 + \binom{n-1}{1}(\dfrac{\beta}{\beta+\gamma})^{n-2} b_2 + \cdots + \binom{n-1}{n-1}(\dfrac{\beta}{\beta+\gamma})^{0} b_n ,\\
g'^{(1)} &=&  b_1.
\end{eqnarray*}
}
Setting
\begin{equation*}
B_{p}=
\begin{pmatrix}
\binom{n-p}{0}      && \binom{n-p}{-1}      && \cdots && \binom{n-p}{-p+1} \\
\binom{n-p}{1}      && \binom{n-p}{0}        && \cdots && \binom{n-p}{-p+2}	\\
\vdots &&  \vdots && \vdots                 && \vdots   			\\
\binom{n-p}{n-p}   && \binom{n-p}{n-p-1} && \cdots && \binom{n-p}{0}		\\
\binom{n-p}{n-p+1}&& \binom{n-p}{n-p}    && \cdots && \binom{n-p}{1}		\\
\vdots                 && \vdots                 && \vdots  && \vdots 			\\
\binom{n-p}{n-1}  && \binom{n-p}{n-2}    && \cdots && \binom{n-p}{n-p}	\\
\end{pmatrix},
\end{equation*}
we compute $Y(p), Y'(P)$ which we want as follows: for $p \neq n-1$,
\begin{eqnarray*}
Y(p) &=& e^{(p)} \wedge f^{(n-p-1)} \wedge g^{(1)} \\
&=& 
\begin{vmatrix}
0 & & \binom{n-1}{0}(\dfrac{\beta}{\beta+\gamma})^{n-1} \\
& &   \binom{n-1}{1}(\dfrac{\beta}{\beta+\gamma})^{n-2} \\
& \mbox{\huge $B_{n-p-1}$ }   & \vdots \\ 
Id_p &   & \binom{n-1}{n-1}(\dfrac{\beta}{\beta+\gamma})^{0}
\end{vmatrix}\\
&=& (-1)^{(n-p)p}
\begin{vmatrix}
\binom{p+1}{0} & \cdots & \binom{p+1}{-n+p+2} & \binom{n-1}{0} (\dfrac{\beta}{\beta+\gamma})^{n-1}\\
\vdots & \vdots & \vdots \\
\binom{p+1}{n-p-1} & \cdots & \binom{p+1}{1} & \binom{n-1}{n-p-1} (\dfrac{\beta}{\beta+\gamma})^{p}
\end{vmatrix}.
\end{eqnarray*}
and $Y(n-1) = (-1)^{n-1}\binom{n-1}{0}(\dfrac{\beta}{\beta+\gamma})^{n-1}$. Similarly for $p \neq n-1$,
\begin{eqnarray*}
Y'(p) &=& e^{(p)} \wedge f^{(n-p-1)} \wedge g'^{(1)} \\
&=&
\begin{vmatrix}
0 & & 1 \\
& &   0 \\
& \mbox{\huge $B_{n-p-1}$ }   & \vdots \\ 
Id_p &   & 0
\end{vmatrix}\\
&=& (-1)^{np+n+1}
\begin{vmatrix}
\binom{p+1}{1} & \cdots & \binom{p+1}{-n+p+3} 			\\
\vdots &  & \vdots												\\
\binom{p+1}{n-p-1} & \cdots & \binom{p+1}{1} 
\end{vmatrix},
\end{eqnarray*}
and otherwise $Y'(n-1) = (-1)^{n-1}$. 
 
\begin{flushleft}
\fbox{5.Shearing invariant $\sigma_p^{\rho_n}(h_{CA})$}
\end{flushleft}

This computation is also similar to the above computation.
We lift $h_{CA}$ to $\widetilde{h}_{CA}$ so that the starting and terminal points are $\infty$ and $1$, and consider the rectangle whose vertices are $1,  0, \infty$ and $\alpha^2 \beta \gamma +1$ where $(\rho(a))(0)  = \alpha^2 \beta \gamma +1$. 
We set $E=\xi_{\rho_n}(1), F=\xi_{\rho_n}(\infty), G=\xi_{\rho_n}(\alpha^2\beta\gamma +1)$ and $G'=\xi_{\rho_n}(0)$. 
Then these flags are given by 
{\small
\begin{eqnarray*}
E^{(p)} &=& {\rm Span}_{\mathbb{R}} <(X+Y)^{n-p}X^{p-1}, (X+Y)^{n-p}X^{p-2}Y, \cdots, (X+Y)^{n-p}Y^{p-1}> ,\\
F^{(p)} &=& {\rm Span}_{\mathbb{R}} <X^{n-1}, X^{n-2}Y, \cdots, X^{n-p}Y^{p-1}> ,\\
G^{(p)} &=& {\rm Span}_{\mathbb{R}} <((\alpha^2 \beta \gamma +1) X+Y)^{n-p}X^{p-1}, ((\alpha^2 \beta \gamma +1) X+Y)^{n-p}X^{p-2}Y, \cdots \\
&& \hspace{10pt} \cdots, ((\alpha^2 \beta \gamma +1) X+Y)^{n-p}Y^{p-1}> ,\\
G'^{(p)} &=& {\rm Span}_{\mathbb{R}} <Y^{n-1}, XY^{n-2}, \cdots, X^{p-1}Y^{n-p}>.
\end{eqnarray*}
}
We choose nonzero elements as
{\small
\begin{eqnarray*}
e^{(p)} &=&  (X+Y)^{n-p}X^{p-1} \wedge (X+Y)^{n-p}X^{p-2}Y \wedge \cdots \wedge (X+Y)^{n-p}Y^{p-1} ,\\
f^{(p)} &=& X^{n-1} \wedge X^{n-2}Y \wedge \cdots \wedge X^{n-p}Y^{p-1}, \\
g^{(p)} &=& ((\alpha^2 \beta \gamma +1) X+Y)^{n-p}X^{p-1} \wedge ((\alpha^2 \beta \gamma +1) X+Y)^{n-p}X^{p-2}Y \wedge \cdots \\
&& \hspace{10pt} \cdots \wedge ((\alpha^2 \beta \gamma +1) X+Y)^{n-p}Y^{p-1},\\
g'^{(p)} &=& Y^{n-1} \wedge XY^{n-2} \wedge \cdots \wedge X^{p-1}Y^{n-p} .
\end{eqnarray*}
}
The vectors $g^{(1)}$ and $g'^{(1)}$ are presented by their entries as follows:
{\small
\begin{eqnarray*}
g^{(1)} &=&  \binom{n-1}{0}(\alpha^2 \beta \gamma +1)^{n-1} b_1 + \binom{n-1}{1}(\alpha^2 \beta \gamma +1)^{n-2} b_2 + \cdots + \binom{n-1}{n-1}(\alpha^2 \beta \gamma +1)^{0} b_n, \\
g'^{(1)} &=&  b_n.
\end{eqnarray*}
}
Then we compute $Y(P)$ and $Y'(P)$ as follows: for  $p \neq 0$

\begin{eqnarray*}
Y(p) &=& e^{(p)} \wedge f^{(n-p-1)} \wedge g^{(1)} \\
&=&
\begin{vmatrix}
& Id_{n-p-1} & \binom{n-1}{0}(\alpha^2 \beta \gamma +1)^{n-1} \\
& &   \binom{n-1}{1}(\alpha^2 \beta \gamma +1)^{n-2} \\
\mbox{\huge $B_{p}$}  & & \vdots \\ 
 &  0 & \binom{n-1}{n-1}(\alpha^2 \beta \gamma +1)^{0}
\end{vmatrix}\\
&=& (-1)^{np}
\begin{vmatrix}
\binom{n-p}{n-p-1} & \cdots & \binom{n-p}{n-2p} &	\binom{n-1}{n-p-1} (\alpha^2 \beta \gamma +1)^p					\\
\vdots &  & \vdots &  \vdots														\\
\binom{n-p}{n-1} & \cdots & \binom{n-p}{n-p} & \binom{n-1}{n-1}  (\alpha^2 \beta \gamma +1) ^0
\end{vmatrix}
\end{eqnarray*}
and $Y(0)=1$. Similarly, for $p \neq 0$
\begin{eqnarray*}
Y'(p) &=& e^{(p)} \wedge f^{(n-p-1)} \wedge g'^{(1)} \\
&=& 
\begin{vmatrix}
 &  Id_{n-p-1} &  0\\
 \mbox{\huge $B_{p}$ } & &    \\ 
& \mbox{\LARGE 0} & 1 
\end{vmatrix}\\
&=& (-1)^{np}
\begin{vmatrix}
\binom{n-p}{n-p-1} & \cdots & \binom{n-p}{n-2p} \\
\vdots & \vdots & \vdots \\
\binom{n-p}{n-2} & \cdots & \binom{n-p}{n-p-1} \\
\end{vmatrix}
\end{eqnarray*}
and $Y'(0) = 1$.

We summarize the result of the computation.

\begin{theorem}
The Bonahon-Dreyer coordinate of the $n$-Fuchsian representations $[\rho_n]$ in the Hitchin component ${\rm Hit}_n(P)$ is explicitly described by using $\alpha, \beta, \gamma$ which are the parameters of the Fuchsian representations in Proposition 5.1 as follows.
\begin{itemize}
\item The shearing invariant $\sigma_p^{\rho_n}(h_{AB})$$(1 \leq p \leq n-1)$
\[ \sigma_p^{\rho_n}(h_{AB}) = \log - \dfrac{Y_{h_{AB}}(p)}{Y'_{h_{AB}}(p)} \cdot \dfrac{Y'_{h_{AB}}(p-1)}{Y_{h_{AB}}(p-1)} \]
where
\[ Y_{h_{AB}}(p) = \binom{n-1}{p}(\beta \gamma)^{n-p-1}, \]
\[ Y'_{h_{AB}}(p) = (-1)^{n-p-1}\binom{n-1}{p} .\]

\item The shearing invariant $\sigma_p^{\rho_n}(h_{BC})$$(1 \leq p \leq n-1)$
\[ \sigma_p^{\rho_n}(h_{BC}) = \log - \dfrac{Y_{h_{BC}}(p)}{Y'_{h_{BC}}(p)} \cdot \dfrac{Y'_{h_{BC}}(p-1)}{Y_{h_{BC}}(p-1)} \] 
where
\[ Y_{h_{BC}}(p) = (-1)^{(n-p)p}
\begin{vmatrix}
\binom{p+1}{0} & \cdots & \binom{p+1}{-n+p+2} & \binom{n-1}{0} (\dfrac{\beta}{\beta+\gamma})^{n-1}\\
\vdots & \vdots & \vdots \\
\binom{p+1}{n-p-1} & \cdots & \binom{p+1}{1} & \binom{n-1}{n-p-1} (\dfrac{\beta}{\beta+\gamma})^{p}
\end{vmatrix}\]
if $p \neq n-1$ and $Y_{h_{BC}}(n-1) = (-1)^{n-1}\binom{n-1}{0}(\dfrac{\beta}{\beta+\gamma})^{n-1}$,
\[ Y'_{h_{BC}}(p) = (-1)^{np+n+1}
\begin{vmatrix}
\binom{p+1}{1} & \cdots & \binom{p+1}{-n+p+3} 			\\
\vdots &  & \vdots												\\
\binom{p+1}{n-p-1} & \cdots & \binom{p+1}{1} 
\end{vmatrix} \]
if $p \neq n-1$ and $Y'_{h_{BC}}(n-1) = (-1)^{n-1}$.

\item The shearing invariant $\sigma_p^{\rho_n}(h_{CA})$$(1 \leq p \leq n-1)$
\[ \sigma_p^{\rho_n}(h_{CA}) = \log - \dfrac{Y_{h_{CA}}(p)}{Y'_{h_{CA}}(p)} \cdot \dfrac{Y'_{h_{CA}}(p-1)}{Y_{h_{CA}}(p-1)} \]
where
\[ Y_{h_{CA}}(p) = (-1)^{np}
\begin{vmatrix}
\binom{n-p}{n-p-1} & \cdots & \binom{n-p}{n-2p} &	\binom{n-1}{n-p-1} (\alpha^2 \beta \gamma +1)^p					\\
\vdots &  & \vdots &  \vdots													\\
\binom{n-p}{n-1} & \cdots & \binom{n-p}{n-p} & \binom{n-1}{n-1}  (\alpha^2 \beta \gamma +1) ^0
\end{vmatrix} \]
if $p \neq 0$ and $Y_{h_{CA}}(0)=1$,
\[ Y'_{h_{CA}}(p) = (-1)^{np}
\begin{vmatrix}
\binom{n-p}{n-p-1} & \cdots & \binom{n-p}{n-2p} \\
\vdots & \vdots & \vdots \\
\binom{n-p}{n-2} & \cdots & \binom{n-p}{n-p-1} \\
\end{vmatrix} \]
if $p \neq 0$ and $Y'_{h_{CA}}(0) = 1$.

\item The triangle invariant $\tau_{pqr}^{\rho_n}(T_0, \infty)$$(p,q,r \geq 1~s. t.~p+q+r=n)$
\[ \tau_{pqr}^{\rho_n}(\widetilde{T}_0, \infty) = \log \dfrac{X_{T_0}(p+1, q, r-1)}{X_{T_0}(p-1, q, r+1)} \cdot \dfrac{X_{T_0}(p, q-1, r+1)}{X_{T_0}(p, q+1, r-1)} \cdot \dfrac{X_{T_0}(p-1, q+1, r)}{X_{T_0}(p+1, q-1, r)} \]
where
\[ X_{T_0}(p, q, r) =
\begin{vmatrix}
\binom{p+r}{p} & \cdots & \binom{p+r}{p-q+1} \\
\vdots & \vdots & \vdots \\
\binom{p+r}{p+q-1} & \cdots & \binom{p+r}{p} \\
\end{vmatrix} \]
if $q \neq 0$ and $X_{T_0}(p,0,r) = 1$ for all $p,r$.

\item The triangle invariant $\tau_{pqr}^{\rho_n}(T_1, \infty)$$(p,q,r \geq 1~s. t.~p+q+r=n)$
\[ \tau_{pqr}^{\rho_n}(\widetilde{T}_1, \infty) = \log \dfrac{X_{T_1}(p+1, q, r-1)}{X_{T_1}(p-1, q, r+1)} \cdot \dfrac{X_{T_1}(p, q-1, r+1)}{X_{T_1}(p, q+1, r-1)} \cdot \dfrac{X_{T_1}(p-1, q+1, r)}{X_{T_1}(p+1, q-1, r)} \] 
where
\[  X_{T_1}(p, q, r) = (-1)^{q(r+1)}
\begin{vmatrix}
\binom{p+q}{p}(-\beta  \gamma)^{q} & \cdots & \binom{p+q}{p-r+1}(-\beta  \gamma)^{q+r-1} \\
\vdots & \vdots & \vdots \\
\binom{p+q}{p+r-1}(-\beta  \gamma)^{q-r+1} & \cdots & \binom{p+q}{p}(-\beta  \gamma)^{q} \\
\end{vmatrix} \]
if $r \neq 0$ and $ X_{T_1}(p,q,0) = (-1)^q$ for all $p,q$. 
\end{itemize}
\end{theorem}
\begin{flushright}
$\Box$
\end{flushright}

\subsection{Remarks on the triangle invariants.}
We give some remarks on the triangle invariants for $n$-Fuchsian representations. 
Recall that the flag curve associated to $n$-Fuchsian representation $\rho_n = \iota_n \circ \rho$ is the composition of the developing map ${\bf dev}_{\rho}$ corresponding to the Fuchsian representation $\rho$ with the map $\xi : \partial_{\infty}\mathbb{H}^2 \rightarrow {\rm Flag}(\mathbb{R}^n)$ defined in the beginning of Section 5.2. 
By definition of the irreducible representation $\iota_n$ and the map $\xi$, it follows that $\iota_n(A) \cdot \xi = \xi \circ A$ for any $A \in {\rm PSL}_2(\mathbb{R})$.
Since the triple ratio is invariant under the action of projective automorphism, we obtain the relation
\begin{eqnarray*}
T_{pqr}(\xi \circ A(x_1), \xi \circ A(x_2), \xi \circ A(x_3)) 
&=& T_{pqr}(\iota_n(A) \xi(x_1), \iota_n(A) \xi(x_2),\iota_n(A) \xi(x_3)) \\
&=& T_{pqr}(\xi(x_1), \xi(x_2), \xi(x_3)) 
\end{eqnarray*}
for any $A \in {\rm PSL}_2(\mathbb{R})$ and any triple $(x_1, x_2, x_3)$ of points in $\partial_{\infty}\mathbb{H}^2$.
Note that the diagonal action of ${\rm PSL}_2(\mathbb{R})$ on the ordered triple $(\partial_{\infty}\mathbb{H}^2)^{(3)}$, the set of triples of pairwise-distinct points of $\partial_{\infty}\mathbb{H}^2$ in clockwise ordering, is transitive.
In particular for any triple $(x_1, x_2, x_3) \in  (\partial_{\infty}\mathbb{H}^2)^{(3)}$ there exists $A \in {\rm PSL}_2(\mathbb{R})$ so that $A \cdot (x_1, x_2, x_3) = (\infty, 1, 0)$.
\begin{proposition}
For any triple $(x_1, x_2, x_3) \in  (\partial_{\infty}\mathbb{H}^2)^{(3)}$, $T_{pqr}(\xi(x_1), \xi(x_2), \xi(x_3)) = T_{pqr}(\xi(\infty), \xi(1), \xi(0))$.
\end{proposition} 
\begin{corollary}
The triangle invariant $\tau_{pqr}^{\rho_n}(T,v)$ of $[\rho_n] \in {\rm Fuch}_n(S)$ always equals to $\log T_{pqr}(\xi(\infty), \xi(1), \xi(0))$, and this value agrees with the following one which is computed in Theorem 5.3.
\[ \log \dfrac{X_{T_0}(p+1, q, r-1)}{X_{T_0}(p-1, q, r+1)} \cdot \dfrac{X_{T_0}(p, q-1, r+1)}{X_{T_0}(p, q+1, r-1)} \cdot \dfrac{X_{T_0}(p-1, q+1, r)}{X_{T_0}(p+1, q-1, r)}. \]
\end{corollary}
{\it Proof.}
We denote the vertices of the ideal triangle $T$ by $v, v'$ and $v''$ which are in clockwise ordering. 
Let $[\rho_n]=[\iota_n \circ \rho] \in {\rm Fuch}_n(S)$ and ${\bf dev}_{\rho}$ be the associated developing map.
Since $({\bf dev}_{\rho}(v), {\bf dev}_{\rho}(v'), {\bf dev}_{\rho}(v'')) \in  (\partial_{\infty}\mathbb{H}^2)^{(3)}$, 
\begin{eqnarray*}
\tau_{pqr}^{\rho_n}(T,v) &=& \log T_{pqr}(\xi({\bf dev}_{\rho}(v)), \xi({\bf dev}_{\rho}(v')), \xi({\bf dev}_{\rho}(v''))) \\
&=& \log T_{pqr}(\xi(\infty), \xi(1), \xi(0)),
\end{eqnarray*}
and we computed this value in theorem 5.3.
\begin{flushright}
$\Box$
\end{flushright}


\begin{thebibliography}{9}
\bibitem{B}
F. Bonahon {\it Shearing hyperbolic srufaces, bending pleated surfaces and Thurston's symplecic form}, Ann. Fac. Sci. Toulouse Math. (6) \textbf{5}(1996), 233-297

\bibitem{BCFS}
M. Bridgeman, R. Canary, F. Labourie and A. Sambarino, {\it The pressure metric for Anosov representations}, Geom. Funct. Anal. \textbf{25}(2015), no. 4, 1089-1179

\bibitem{BD1}
F. Bonahon and G. Dreyer, {\it Hitchin characters and geodesiic laminations}, arXiv:14100729

\bibitem{BD2}
F. Bonahon and G. Dreyer, \textit{Parameterizing Hitchin components}, Duke Math. J. \textbf{163}(2014), no. 15, 2935-2975.

\bibitem{CG}
S. Choi and W. Goldman, {\it Convex real projective structures on closed surfaces are closed}, Proc. Amer. Math. Soc. \textbf{118}(1993), no. 2, 657-661.

\bibitem{FG}
V. Fock and A. Goncharov, \textit{Moduli spaces of local systems and higher Teichm\"uller theory}, Publ. Math. Inst. Hautes \'Etudes Sci. No. \textbf{103}(2006) 1-211.

\bibitem{Go}
W. Goldman, {\it Topological components of spaces of representations}, Invent. Math. \textbf{93}(1988), no. 3, 557-607.

\bibitem{Gu}
O. Guichard, \textit{Composantes de Hitchin et repr\'esentations hypervonvexes de groupes de surface}, J. Differential Geom. \textbf{80}(2008), no. 3, 391-431.

\bibitem{GGKW}
F. Gu\'eritaud, O. Guichard, F. Kassel and A. Wienhard, {\it Anosov representations and proper actions}, Geom. Topol. \textbf{21}(2017), no. 1, 485-584.

\bibitem{GW1}
O. Guichard and A. Wienhard, {\it Anosov representations: domains of discontinuity and applications}, Invent. Math. \textbf{190}(2012), no. 2, 357-438.

\bibitem{GW2}
O. Guichard and A. Wienhard, {\it Convex foliated projective structures and the Hitchin component for ${\rm PSL}_4(\mathbb{R})$}, Duke Math. J. \textbf{144}(2008), no. 3, 381-445

\bibitem{Hi}
N. Hitchin, {\it Lie groups and Teich\"uller space}, Topology \textbf{31}(1992), no. 3, 449-473  

\bibitem{IT}
Y. Imayoshi and M. Taniguchi, {\it An introduction to Teichm\"uller spaces}, Springer-Verlag, Tokyo, 1992.

\bibitem{KLP}
M. Kapovich, B. Leeb and J. Porti, {\it Morse actions of discrete groups on symmetric space}, arXiv:14037671

\bibitem{L}
F. Labourie, \textit{Anosov flows, surface groups and curves in projective space}, Invent. Math. \textbf{165}(2006), no. 1, 51-114.

\bibitem{LM}
F. Labourie and G. McShane, \textit{Cross ratios and identities for higher Teichm\"uller-Thurston theory}, Duke Math. J. \textbf{149}(2009), no. 2, 279-345.

\bibitem{LW}
F. Labourie and R. Wentworth, {\it Variations along the Fuchsian Locus}, arXiv:150601686

\bibitem{PS}
R. Potrie and A. Sambarino, {\it Eigenvalues and entropy of a Hitchin representation}, Invent. Math. \textbf{209}(2017), no. 3, 885-925.

\bibitem{R}
J. Ratcliffe, {\it Foundations of hyperbolic manifolds}, Graduate Texts in Math. \textbf{149}, Springer-Verlag, New York, (2006).

\bibitem{S}
A. Sikora, {\it Character varieties}, Trans. Amer. Math. Soc. {\bf 364}(2012), no. 10, 5173-5208.

\bibitem{T}
W. Thurston, {\it Minimal stretch maps between hyperbolic surfaces}, arXiv:9801039

\bibitem{Z}
T. Zhang, {\it Degeneration of Hitchin representations along internal sequences}, Geom. Funct. Anal. \textbf{25}(2015), no. 5, 1588-1645.

\end{thebibliography}
\end{document}